\documentclass[reqno,11pt]{amsart}
\usepackage{tkz-euclide}
\usepackage{hyperref}
\usepackage{amsmath} 
\usepackage{amsfonts}
\usepackage{mathtools}
\usepackage{amssymb} 
\usepackage{mathrsfs} 
\usepackage{graphicx}  
\usepackage{cancel}
\usepackage{bm}
\usepackage{tikz} 
\usepackage{centernot}
\usepackage{soul}
\usetikzlibrary{svg.path} 
\usepackage{longtable}
\usepackage{makecell}

\usepackage{geometry}
\usetikzlibrary{calc} 
\usepackage{xcolor}  
\usetikzlibrary{patterns}
\usetikzlibrary{arrows.meta} 
\usepackage{amsthm} 
\usepackage{float}
\usepackage{enumitem}
\setenumerate{ ref={\normalfont(\roman*)},label={\normalfont(\roman*)},itemsep=3pt}
\usepackage[english]{babel}
\usepackage[font=footnotesize, labelfont=bf]{ caption}
\usepackage{ esint }
\usepackage{stackengine,scalerel,graphicx}
\stackMath

\definecolor{trueblue}{rgb}{0.0, 0.45, 0.81}
\definecolor{truegreen}{rgb}{0.13, 0.55, 0.13}

\theoremstyle{plain}
\begingroup
\newtheorem{theorem}{Theorem}[section]
\newtheorem{lemma}[theorem]{Lemma}

\newtheorem{proposition}[theorem]{Proposition}
\newtheorem{corollary}[theorem]{Corollary}

\endgroup

\theoremstyle{definition}
\begingroup
\newtheorem{definition}[theorem]{Definition}
\newtheorem{remark}[theorem]{Remark}
\endgroup

\renewcommand{\tilde}{\widetilde}

\numberwithin{equation}{section}
\newcommand{\N}{\mathbb{N}}
\newcommand{\Z}{\mathbb{Z}}

\newcommand{\E}{\mathcal{E}}
\newcommand{\R}{\mathbb{R}}
\newcommand{\C}{\mathbb{C}}

\newcommand{\I}{\mathcal{I}}

\newcommand{\BBB}{\color{black}}

\newcommand{\EEE}{\color{black}} 
\newcommand{\PPP}{\color{black}}

\renewcommand{\S}{{\mathbb{S}}}

\renewcommand{\H}{\mathcal{H}}

\usepackage[normalem]{ulem}

\def \mb{\mathbb}
\def \Id{\textup{Id}}

\def \mc{\mathcal}

\def \tp{\textup}

\begin{document}
	
\title
[ Global minimality  of the Hopf map in the Faddeev--Skyrme model]
{ Global minimality  of the Hopf map in the Faddeev--Skyrme model with large coupling constant}
	
\author[Andr\'e Guerra]{Andr\'e Guerra} 
\address[Andr\'e Guerra]{Department of Pure Mathematics and Mathematical Statistics, University of Cambridge, Wilberforce Rd, Cambridge CB3 0WB, UK}
\email{adblg2@cam.ac.uk}

\author[Xavier Lamy]{Xavier Lamy} 
\address[{Xavier Lamy}]{Institut de Mathématiques de Toulouse, Université de Toulouse, 118 route de Narbonne, F-31062 Toulouse Cedex 8, France}
\email{Xavier.Lamy@math.univ-toulouse.fr}	
	
\author{Konstantinos Zemas}
\address[Konstantinos Zemas]{Institute for Applied Mathematics, University of Bonn\\
Endenicher Allee 60, 53115 Bonn, Germany}
\email{zemas@iam.uni-bonn.de}
	

\begin{abstract}
We prove that, modulo rigid motions, the Hopf map is the unique minimizer of the Faddeev--Skyrme energy in its homotopy class,
 \BBB provided that the radius of the target \PPP 2-\BBB sphere is not smaller than the radius of the domain \PPP 3-\BBB sphere.  \EEE
\end{abstract} 
	
\maketitle
\thispagestyle{empty}


\section{Introduction}

The Skyrme model, introduced in \cite{skyrme1961non},  is a classical nonlinear $O(3)$-sigma model which has proved successful in quantum field theory. The static fields in Skyrme's model, which are called \textit{Skyrmions}, can be mathematically described  as critical points of the energy
\[ \mc S(u) := \int_{\mb{M}^3}|du|^2+\frac{1}{4}\int_{\mb{M}^3}|du\wedge du|^2\,,\qquad u\colon \mb M^3\to \mb S^3\,,\]
where $\mb M^3$ is a Riemannian manifold\footnote{Here, $\wedge$ is the natural wedge 
product on $T^*\mb M^3\otimes \R^4$ characterized by $(\alpha\otimes a) \wedge (\beta\otimes b )=(\alpha\wedge\beta)\otimes (a\wedge b)$.}.
Skyrmions are a particular example of topological solitons, and we refer the reader to the monograph \cite{Manton2004} for further physical background. The existence of maps minimizing $\mc S$ in suitable homotopy  classes 
has been proved rigorously in the works \cite{Esteban1986,Esteban2004}, in the case $\mb M^3=\R^3$. 
In this case, finite-energy maps exhibit rapid decay at infinity to a constant value $b\in \S^3$, so that they can be considered as maps from  $\R^3\cup\{\infty\}\cong \S^3\to \S^3$. 

A further refinement of Skyrme's model was proposed by Faddeev \cite{faddeev1976some},  considering 
Skyrme's model in the case where the maps take values in an equator $\mb S^2\subset \mb S^3$. As in Skyrme's model, there is a satisfactory existence theory for this problem, especially in the case of compact $\mb{M}^3$ \cite{Auckly2005,lin2004existence,lin2007energy, ward1999}, the case of non-compact $\mb{M}^3$ being more challenging \cite{lin2004existence,lin2007energy}.  In the simplest non-compact example $\mb M^3=\R^3$, and following Ward \cite{ward1999},  one  can consider an approximation of $\R^3$ by $\mb S^3_\rho$,  the sphere of radius $\rho\gg 1$. Thus, considering the energy on $\S^3_\rho$, after a change of variables 
 and normalization by a factor $\rho^{-1}$, one is led to a \textit{perturbed version} of the Faddeev--Skyrme energy,  namely  
\begin{equation}\label{eq:pert_Faddeev_skyrme}
\mathcal{FS}_\rho(u):= \int_{\mathbb S^3}|du|^2 + \frac{1}{4\rho^{2}}\int_{\mathbb S^3}|du \wedge d u|^2\,, \qquad u\colon \mb S^3\to \mb S^2.
\end{equation}
One can rewrite the second term in the energy above as 
$$\frac 1 4 |du\wedge d u|^2=|u^*\omega_{\mb S^2}|^2\,,$$ 
where $\omega_{\mb S^2}$ is the volume form  on  
$\mb S^2$, 
 see 
\eqref{eq:modulus_ofu_omega} below. 
The energy \eqref{eq:pert_Faddeev_skyrme} is therefore naturally defined in the space
\begin{equation}\label{eq:energy_space}
\mathcal{U}_{\mathcal{FS}}:=\left\{u\in W^{1,2}(\S^3;\S^2)\colon \int_{\S^3}|du|^2+\int_{\S^3}|u^*\omega_{\S^2}|^2<+\infty\right\}\,.
\end{equation} 
Homotopy classes of smooth maps from $\mb S^3$ to $\S^2$,
corresponding to elements of $\pi_3(\mathbb S^2)\cong \mb Z$ are classified by their 
\textit{Hopf invariant} 
 (or \textit{Hopf charge} in the physics literature) \cite{whitehead1947expression}, which 
is  defined for $u\in C^\infty(\S^3;\S^2)$ as 
\begin{equation}\label{eq:1_Hopf_inv}
Q(u):=\frac 1{16\pi ^2} \int_{\mathbb S^3} \beta\wedge d\beta\in \Z\,, 
\end{equation}
where $\beta$ is 
any 1-form on $\S^3$ such that $d\beta=u^\ast\omega_{\S^2}$. This definition can be meaningfully extended to the class of finite Faddeev--Skyrme energy maps, see \cite{auckly2010pontrjagin} or Section \ref{sec:degree} below.

A basic question concerning the above models is whether one can characterize minimizers in some simple cases; the non-convexity of the problems makes this a rather intriguing question. For the Skyrme model with $\mb M^3=\mb S^3_\rho$, it is conjectured \cite{Manton1986} that the homothety $x\mapsto x/\rho$ is a minimizer in its homotopy class whenever $ 0 <\rho \leq \sqrt 2$; this map is known to be unstable whenever $\rho > \sqrt 2$. At present minimality and uniqueness modulo isometries is only known for 
$ 0<\rho \leq \sqrt{3/2}<\sqrt{2}$, see \cite{riviere1998remark} and the references therein.

For the Faddeev--Skyrme model \eqref{eq:pert_Faddeev_skyrme} the situation is more complicated, as one has \textit{knotted} solutions to the variational equations, and the \textit{Hopf map} $h$ is a critical point of $\mathcal{FS}_\rho$ of particular interest from the point of view of both geometry and  physics. This map gives a prominent example of a non-trivial $\mb S^1$-fiber bundle of $\S^3$ onto $\S^2$: in complex coordinates, and denoting by $\pi\colon \C \to \S^2$ the inverse stereographic projection, \textit{i.e.},
\begin{equation*}
\pi\colon \C\cup\{\infty\} \to \S^2\,, \quad  \pi(z):=\left(\frac{2z}{|z|^2+1},\frac{|z|^2-1}{|z|^2+1}\right)\,,
\end{equation*}
the Hopf map $h:\S^3\subset \C^2\to \S^2$ is given by
\begin{equation}\label{def:Hopf}
h(z,w):=\pi(z/w)=(2z\bar w, |z|^2-|w|^2)\,,
\end{equation}
where the last equality holds true because $|z|^2+|w|^2=1$ for $(z,w)\in \S^3$.  

For the Faddeev--Skyrme model it is also known \cite{ward1999}  that $h$ is unstable  for $\rho>\sqrt{2}$, while it is linearly stable for $ 0<\rho\leq \sqrt{2}$,  cf. \cite{isobe2008,speight2007}, similarly to the homothety in the Skyrme model. It is therefore conjectured that, in the latter case, $h$ is the unique (modulo isometries) minimizer of the energy in its homotopy class.

Actually, as pointed out\footnote{Specifically, cf.\ the paragraph below Theorem 1.4 in \cite{speight2011}.} in \cite{speight2011}, the above conjectured minimality property of the Hopf map  was unknown  
for \textit{any} value of the coupling constant 
$ 0<\rho\leq \sqrt 2$. In this paper we give a positive answer to the above conjecture \BBB for a certain range of $\rho$: 
\begin{theorem}\label{thm:main}
\BBB For every $\rho\in (0,\PPP 1 ]\EEE$ and $u\in \mathcal{U}_{\mathcal{FS}}$ with Hopf invariant equal to 1, it holds that
\begin{equation*}\label{eq:minimality_hopf}
\mathcal{FS}_\rho(u) \geq \mathcal{FS}_\rho(h)\,,
\end{equation*}
with equality if and only if $u=h\circ R$ for some $R\in \textup{SO}(4)$.\EEE
\end{theorem}
\EEE

The case where the energy consists only of the second term in \eqref{eq:pert_Faddeev_skyrme}, which formally corresponds to the limit $\rho\to 0$, was previously considered in \cite{speight2011}, where a short proof of the minimality of the Hopf map was provided.

It is also interesting to compare our results with another natural conjecture, tha the Hopf map should minimize, in its homotopy class, the conformally-invariant energy
\begin{equation}
\label{eq:riv}
\int_{\mb S^3} |du|^3\,, \qquad u\colon \mb S^3\to \mb S^2\,.
\end{equation}
Linear stability of $h$  in this case  was already proved in \cite{riviere1998}, as well as minimality if the power 3 is replaced with 4, but otherwise this conjecture remains open; see also \cite{CoronGulliver} for related results. On the one hand, there are several similarities between \eqref{eq:riv} and the Faddeev--Skyrme problem,  cf.\ for instance  \cite{lin2004existence,lin2007energy}, and as we will see in this paper the required spectral analysis is also quite similar. On the other hand, the problems are at the same time fairly different, since \eqref{eq:pert_Faddeev_skyrme} is not conformally invariant and the corresponding integrand is non-convex.

\PPP A natural strategy to prove Theorem~\ref{thm:main} 
would elaborate on two already mentioned known facts about the Hopf map:
\begin{itemize}
\item[(a)] it is linearly stable for $0<\rho\leq \sqrt 2$, cf.  \cite{isobe2008,speight2007};
\smallskip
\item[(b)] it is minimal for the second term in the Faddeev-Skyrme energy \eqref{eq:pert_Faddeev_skyrme},  cf.\  \cite{speight2011}.
\end{itemize}
If one could establish quantitative versions of these two  facts: loosely stated that,
\begin{itemize}
\item[(a')] the Hopf map is (modulo symmetries) a strict $W^{1,2}$-local minimizer for
$0<\rho<\sqrt 2$;
\smallskip 
\item[(b')] if a map is almost minimal for the second energy term, then it must be close to the Hopf map (modulo symmetries),
\end{itemize}
then since the second energy term in \eqref{eq:pert_Faddeev_skyrme} is dominant for small $\rho$,
these strengthened  assertions would provide a proof of Theorem~\ref{thm:main} for a sufficiently small, though non-explicit range of $\rho$. Our proof is actually a more direct variant of the aforementioned strategy and it yields an explicit range. A key difficulty is related to the symmetries in these statements: the second energy term in \eqref{eq:pert_Faddeev_skyrme}, which
depends only on the closed 2-form $\alpha=u^*\omega_{\S^2}$, has more symmetries than the full energy. 
But this difficulty is resolved thanks to
an idea introduced by Rivière \cite{riviere1998,riviere1998remark} and used also in \cite{isobe2008},
which allows to apply the outlined strategy 
using a relaxed energy on closed $2$-forms. Note that, in contrast to our setting, in \cite{riviere1998remark} the relaxed energy is quadratic, cf.\ (1.10) therein, and thus the associated \BBB relaxed \PPP minimization problem can be solved explicitely. \EEE

{\bf Plan of the paper:} In Section \ref{sec:prelim} we gather some notation and basic properties related to the Hopf map and the Hodge Laplacian on $\mb S^3$. In Section \ref{subsec:Hopf_inv} we prove several useful facts about the Hopf invariant for maps from $\S^3$ to $\S^2$ with finite Faddeev--Skyrme energy, and 
 include  
a concise proof of its integrality in this lower-regularity case. In Section \ref{sec:relaxed_formulation}, \PPP we establish 
global minimality properties for a relaxed energy functional on closed $2$-forms, \PPP the main result therein being Proposition \ref{prop:stability}. With this as our key ingredient, \PPP at the end of the section we obtain the minimality of the Hopf map as an immediate consequence, and in \EEE Section \ref{sec:proof_main} we conclude the proof of Theorem \ref{thm:main} \PPP by verifying the uniqueness (modulo rotational symmetries) in the statement\EEE.

\section{Preliminaries}\label{sec:prelim}
For $p\in [1,+\infty)$ and for a vector bundle $E\to \S^3$, we denote by $L^p(\S^3;E)$ (respectively $W^{1,p}(\S^3;E)$) the space of $L^p$- (respectively $W^{1,p}$-) sections of $E\to \S^3$. 

\subsection{Some basic identities involving differential forms}
\label{subsec:diff_forms}

Let us first define precisely the terms in \eqref{eq:pert_Faddeev_skyrme}. Here, we adopt an \textit{extrinsic viewpoint}, namely we consider $$u:=(u^1,u^2,u^3):\S^3\to \mb S^2 \subset \R^3\,,$$ 
where $\R^3$ is endowed with the standard Euclidean basis $(e_1,e_2,e_3)$, so that $du(x)\in T_x^\ast\S^3\otimes \R^3$ is given in coordinates as
\[du(x)=\sum_{k=1}^3du^k(x)\otimes e_k
{=
\begin{pmatrix}du^1\\[2pt]
du^2\\[2pt]
du^3
\end{pmatrix}
}
\,.\] 
This gives rise to an alternating vector-valued 2-form $du(x)\wedge du(x)\in {\bigwedge^2 \PPP(\EEE T_x^\ast\S^3\otimes \R^3\PPP)\EEE}$, which in coordinates is expressed as
\begin{equation}\label{eq:du_wedge_prod}
du(x)\wedge du(x)=\left(\sum_{k=1}^3du^k(x)\otimes e_k\right)\wedge \left(\sum_{k'=1}^3du^{k'}(x)\otimes e_{k'}\right)=
2\begin{pmatrix}du^2\wedge du^3\\[2pt]
du^3\wedge du^1\\[2pt]
du^1\wedge du^2
\end{pmatrix}\,,
\end{equation}
where we also used the relations $e_1\wedge e_2=e_3\,,\ e_2\wedge e_3=e_1\,,\ e_3\wedge e_1=e_2\,.$

For later reference, we note that the volume form on $\mb S^n\subset\R^{n+1}$ is given by 
\begin{equation}
\label{eq:omega_coordinate_expression}
\omega_{\mb S^n} = \sum_{j=1}^{n+1} (-1)^{j-1} x^j dx^1\wedge \dots \wedge \widehat{dx^j} \wedge \dots \wedge dx^{n+1}\,,
\end{equation}
where as usual $\widehat{\cdot}$ denotes that the corresponding term is omitted.
Thus, when $n=2$, its pull-back by $u$ 
 is 
\begin{align}\label{eq:omega_coordinates}
\begin{split}	
u^*\omega_{\S^2}&=u^1du^2\wedge du^3+u^2du^3\wedge du^1+u^3du^1\wedge du^2\\
&=\frac{1}{2}u\cdot du\wedge du=\sum_{1\leq j<k\leq 3}u\cdot(\partial_j u\times \partial_k u)\,dx^j\wedge dx^k\,,	
\end{split}
\end{align}
where $\cdot$ and $\times$ denote the Euclidean inner and outer products  respectively. Since $|u|= 1$, \eqref{eq:du_wedge_prod} and \eqref{eq:omega_coordinates} imply that 
\begin{align}\label{eq:modulus_ofu_omega}
|u^*\omega_{\S^2}|^2=\sum_{1\leq j<k\leq 3}|\partial_j u\times \partial_k u|^2=\frac{1}{4}|du\wedge du|^2\,.	
\end{align}

Throughout this paper we will use systematically Hodge theory for differential forms with coefficients in Sobolev spaces, and we refer the reader to  \cite[\textsection 5]{iwaniec1999nonlinear} for a comprehensive discussion of this theory.  Here we record the following simple lemma.

\begin{lemma}\label{lem:nulllag}  The following statements hold true:
\begin{itemize}
\item[$\rm{(i)}$] Let $\alpha,\beta \in W^{1,2}(\mb S^3;\bigwedge^1 T^*\mb S^3)$. Then,
\begin{equation}\label{eq:forms_int_by_parts}
\int_{\mb S^3} \alpha \wedge d \beta = \int_{\mb S^3}  d \alpha\wedge \beta\,.
\end{equation}
\item[$\rm{(ii)}$] There exists a constant $C>0$ such that for every $\varphi \in L^2(\mb S^3; \bigwedge^1 T^*\mb S^3)$ with $d^*\varphi=0$, there holds
\begin{equation}\label{eq:forms_Poincare}
\Big|\int_{\mb S^3} \varphi \wedge d \varphi\Big| \leq C \int_{\S^3}|d \varphi|^{2}\,.
\end{equation}
\end{itemize}
\end{lemma}

\begin{proof}
The proof of (i) is immediate by the formula $d(\alpha\wedge \beta)=d\alpha\wedge \beta-\alpha\wedge d\beta$ and Stokes' theorem. Regarding (ii), let $\overline\varphi:=\fint_{\S^3} \varphi\in \bigwedge^1 T^*\S^3$ be the \textit{mean value of} $\varphi$, intended component-wise. Then
\begin{align*}
\Big|\int_{\mb S^3} \varphi \wedge d \varphi\Big| &=\Big|\int_{\mb S^3} (\varphi-\overline \varphi) \wedge d \varphi\Big| \leq \|\varphi-\bar \varphi\|_{L^2(\S^3)}\,\|d\varphi\|_{L^2(\S^3)}\\
&\leq C\|\nabla \varphi\|_{L^2(\S^3)}\,\|d\varphi\|_{L^2(\S^3)}\leq C\|d\varphi\|^2_{L^2(\S^3)}\,,
\end{align*}
where in the first line we used that, by \eqref{eq:forms_int_by_parts}, $\int_{\S^3}c\wedge d\varphi=0$ for every constant coefficient 1-form $c\in \bigwedge^1 T^*\S^3$ and the Cauchy-Schwartz inequality. In the passage to the second line we used successively the Poincar{\'e} inequality on forms, that $\varphi$ is co-closed and \textit{Gaffney's inequality}
\begin{align}\label{eq:gaffney}
\|\nabla\varphi\|_{L^2(\S^3)}\leq C\big(\|d\varphi\|_{L^2(\S^3)}+\|d^*\varphi\|_{L^2(\S^3)}\big)\,,
\end{align}
which the reader can find in 
\cite[{\textsection 4}]{iwaniec1999nonlinear}.
\end{proof}

\PPP
Actually, in the proof of Proposition \ref{prop:faddeev}, cf.\ \eqref{eq:psi_dpsi_est} therein, we will calculate the \BBB optimal \PPP value of the constant $C>0$ in \eqref{eq:forms_Poincare}. \EEE
We will also use, often implicitly, the characterization of the De Rham cohomology groups of spheres. In particular, for $\mb S^3$ we have that
\begin{equation*}\label{eq:cohom_groups}
H^k_{\rm{dR}}(\S^3)=\begin{cases}
\R\,, \ \ \text{for }k=0,3\\
0\,, \ \ \text{for }k=1,2\,.
\end{cases}
\end{equation*}
 cf.\
\cite[Theorem~17.21]{lee2003smooth}.
In particular closed 1- and 2-forms on $\mathbb S^3$ are exact,  and the same holds  for forms with Sobolev coefficients, cf.\ \cite[\textsection 5]{iwaniec1999nonlinear}.

\subsection{Spectral properties of the Hodge-Laplacian on closed forms}\label{app:spectrum_hodge}

We collect here some well-known spectral properties of the Hodge-Laplacian $\Delta := d d^*+d^*d$ on closed 2-forms on $\mb S^3$, where we recall that
\[d^* := * d * \colon \Omega^2(\mb S^3)\to \Omega^{ 1}(\mb S^3)\,.\]
In fact, it is also useful to consider the square root of $\Delta$ which, restricted to closed 2-forms, becomes the operator $d*$. The corresponding eigenvalue equation is then
\begin{equation}\label{eq:eigenval}
d * \alpha = \lambda \alpha\,,\ \ \  \text{ with } \alpha \in \Omega^2(\S^3) \text{ such that } d \alpha =0\,.
\end{equation}
Note that if $\alpha$ solves \eqref{eq:eigenval} then $\Delta \alpha = \lambda^2\alpha$. Since the spectrum of $\Delta$ on $\Omega^2(\mb S^3)$ is the set $\{(k+2)^2\}_{k\in \N}$,  cf. \cite{ikeda1978,Iwasaki1979}, it follows that for \eqref{eq:eigenval} to be solvable we must have $\lambda= \pm(k+2)$. In fact \eqref{eq:eigenval} is solvable for all such values,  see again \cite{ikeda1978,Iwasaki1979} as well as \cite[Proposition IV.2]{riviere1998}, and we have the following precise description of the corresponding eigenspaces:
\begin{proposition}[Spectrum of $d\, *$]\label{prop:spectrum}
Let $\Id\colon \Omega^2(\mb S^3)\to \Omega^2(\mb S^3)$ be the identity operator. The eigenspaces
\begin{equation}\label{eq:eigenspaces}
E_k^\pm := \ker\big(d * \mp (k+2)\mathrm{ Id}\big)\cap \ker d\subset \Omega^2(\mb S^3)
\end{equation} 
are the restrictions to $\mb S^3$ of the {self-}dual (respectively anti-self-dual) closed and co-closed 2-forms on $\R^4$ which are homogeneous and polynomial of degree $k$. These spaces induce an $L^2$-orthogonal decomposition
\begin{equation}\label{eq:L_2_orthogonal_Hodge}
L^2\Big(\mb S^3; \bigwedge\nolimits^2 T^* \mb S^3\Big) = \bigoplus_{k=0}^\infty (E_k^+ \oplus E_k^-)\,.
\end{equation}
\end{proposition}

We recall that a form $\alpha\in \Omega^2(\R^4)$ is said to be {self-}dual/anti-self-dual if $* \alpha = \pm \alpha$.

\begin{remark}\label{rem:structure_of_E_0} Since $**=\mathrm{Id}$ on $\Omega^2(\R^4)$ and the Hodge-star operation is an isometry, the spaces of \textit{self-dual} and \textit{anti-self-dual} forms give a pointwise-orthogonal decomposition of $\Omega^2(\R^4)$. In the case $k=0$, the elements in $E_0^\pm$ have a specially simple description, since they are the restrictions to $\mb S^3$ of constant-coefficient 
2-forms on $\R^4$, and hence they can be identified with the self-dual (or anti-self-dual) constant-coefficient 
2-forms on $\R^4$. For instance, we have
\begin{equation}\label{eq:basis_for_E_0+}
E_0^+ = \mathrm{span}\lbrace
\omega^+_{0,1},\ \omega^+_{0,2},\ \omega^+_{0,3}\}\,,
\end{equation}
where
$$\omega^+_{0,1}:=dx^1\wedge dx^2 +dx^3\wedge dx^4\,,\ \omega^+_{0,2}:=dx^1\wedge dx^3  - dx^2\wedge dx^4\,,\ \omega^+_{0,3}:=dx^1\wedge dx^4 +dx^2\wedge dx^3\,,$$
with a similar basis for $E_0^-$, up to flipping the middle signs in the above basis.  Note in particular that $\dim E_0^\pm=3$. 

We will also  need  the fact that the natural action of $\tp{SO}(4)$ on the space of self-dual (or anti-self-dual)  constant-coefficient  2-forms on $\R^4$ is \textit{transitive}, \textit{i.e.},
\begin{equation}\label{SO4_transitive_action}
\forall \omega_1,\omega_2\in E_0^\pm
\text{ with }|\omega_1|=|\omega_2| 
\text{ there exists }
\tilde R\in SO(4) \text{ such that } \omega_2=\tilde R^*\omega_1\,.  
\end{equation}
This follows directly from the fact that any $\omega\in E_{0}^+$ can be represented as $\omega =\lambda R^*\omega_{0,1}^+$ for some $\lambda\in\R$ and $R\in \tp{SO}(4)$ (and analogously for $E_{0}^-$). 
To check this, consider 
the skew-symmetric matrix $A\in \R^{4\times 4}$ of the 2-form $\omega$ in the canonical basis, 
\textit{i.e.}, $A_{ij}:=\omega(e_i,e_j)\ \forall i,j=1,\dots,4$:
there exists a rotation $R\in \tp{SO}(4)$ such that $RAR^t$ is block-diagonal with $2\times 2$ skew-symmetric blocks.
Since $\omega$ is self-dual, 
these two blocks must be the same.
This means that $RAR^t=\lambda A_{0,1}^+$
for some $\lambda\in\R$,
where $A_{0,1}^+$ is the matrix of $\omega_{0,1}^+$, and therefore
 $\omega =\lambda R^*\omega_{0,1}^+$.
\end{remark}

\subsection{The Hopf map and horizontal conformality}\label{subsec:horiz_conf}
In this subsection we gather some well-known facts about the Hopf map.
We first recall the following definition,  which the reader can find more generally in \cite[Definition 2.4.2 and Lemma 2.4.4]{baird2003harmonic}:

\begin{definition}\label{def:horiz_weakly_conf}
A map $u\in W^{1,2}(\S^3;\S^2)$ is said to be \textit{horizontally weakly conformal} if 
\begin{equation}\label{eq:horiz_weak_conf_alternative}
du(x)\circ du(x)^t=\frac{1}{2}|du(x)|^2\,\mathrm{Id}_{T_{u(x)}\S^2} \quad \text{for }\H^3\text{-a.e. } x\in \S^{\PPP 3\EEE}\,,
\end{equation}
where $du(x)^t \colon T_{u(x)}\S^2\to T_x\S^3$ denotes the adjoint map of $du(x)\colon T_x\S^3\to T_{u(x)}\S^2$.
\end{definition}

The following lemma gives a simple but very useful pointwise inequality related to weak horizontal conformality. Its proof is elementary and can be found for instance in \cite[Proof of Lemma 2.1]{isobe2008}.
\begin{lemma}\label{lem:dirichlet_energy_weak_conformality}
For a map $u\in W^{1,2}(\S^3;\S^2)$ the following inequality holds for $\H^3$-a.e. $x\in \S^3\mathrm{:}$
\begin{equation}\label{eq:am_gm-ineq}
|u^*\omega_{\S^2}(x)|\leq \frac{1}{2}|du(x)|^2\,.
\end{equation}
Equality holds if and only if $u$ is horizontally weakly conformal at $x$.
\end{lemma}

We remark that this type of estimate also plays an important role in the minimality arguments in \cite{CoronGulliver}.	

The Hopf map $h$ introduced in \eqref{def:Hopf} is a main example of a horizontally conformal map. We now detail some of its basic properties. For the 1-form on $\S^3\subset\R^4$ given by
\begin{equation}\label{eq:alpha_exp}
\theta:=-x^2dx^1+x^1dx^2-x^4dx^3+x^3dx^4\,,
\end{equation}
it holds that 
\begin{equation}\label{eq:alpha_hopf}
h^*\omega_{\S^2}=2d\theta=4(dx^1\wedge dx^2+dx^3\wedge dx^4)\,,
\end{equation}
and as a result of \eqref{eq:alpha_exp}, \eqref{eq:alpha_hopf} and \eqref{eq:omega_coordinate_expression}  we also have
\begin{equation}\label{eq:a_wedge_da}
\theta \wedge d\theta=2\omega_{\S^3}\,.
\end{equation}

It is actually convenient to choose particular coordinates on $\S^3$ and $\S^2$, namely
\begin{align}\label{eq:local_param}
\begin{split}
\mathrm{(i)}&\quad [0, \pi/2]\times \S^1\times \S^1\ni (t,e^{i\varphi_1}, e^{i\varphi_2})\mapsto (e^{i\varphi_1}\sin t, e^{i\varphi_2}\cos t) \in \S^3\,,\\
\mathrm{(ii)}&\quad [0, \pi/2]\times \S^1\ni (t, e^{i\varphi})\mapsto (e^{i\varphi} \sin 2t,-\cos 2t) \in \S^2\,.
\end{split}
\end{align}
In these coordinates,
\begin{align*}
g_{\S^3}&=(dt)^2 + \sin^2t\,(d\varphi_1)^2+\cos^2t\,(d\varphi_2)^2\,,\\
g_{\S^2}&= 4(dt)^2+(\sin^22t)\,(d\varphi)^2\,,
\end{align*}
and  the Hopf map takes the simple form
\[h\colon (t,e^{i\varphi_1}, e^{i\varphi_2})\mapsto (t, e^{i(\varphi_1-\varphi_2)})\,.
\]
One can choose an orthonormal frame $\{\tau_1, \tau_2, \tau_3\}$, defined at a chart point $(t, e^{i\varphi_1}, e^{i\varphi_2}) \in \S^3$ as
\begin{equation}\label{eq:special_S3_frame}
\tau_1:=\frac{\partial}{\partial \varphi_1}+\frac{\partial}{\partial \varphi_2}\,,\quad \tau_2:=\frac{\partial}{\partial t}\,,\quad \tau_3:=\cot t\frac{\partial}{\partial \varphi_1}-\tan t\frac{\partial}{\partial \varphi_2}\,,
\end{equation}
so that $\tau_1$ is the \textit{fundamental vertical vector ﬁeld}, \textit{i.e.}, the generating vector ﬁeld of the $\S^1$-action $\S^1 \times \S^3\ni (e^{it}, (u,v))\to (e^{it}u,e^{it}v)\in \S^3$. With this choice, 
\begin{equation}\label{eq:dh_coord}
dh(\tau_1)=0\,,\quad dh(\tau_2)=\frac{\partial}{\partial t}\,,\quad dh(\tau_3)=\frac{2}{\sin 2t}\frac{\partial}{\partial \varphi}\,,
\end{equation}
where $\varphi:=\varphi_1-\varphi_2$. In particular, $\{f_1,f_2\}:=\big\{\frac{1}{2}\frac{\partial}{\partial t},\frac{1}{\sin 2t}\frac{\partial}{\partial \varphi}\big\}$ is an orthonormal frame on $\S^2$, and in these coordinates the differential of the Hopf map $dh\colon T\S^3\to T\S^2$ has matrix representation 
\begin{equation}\label{eq:matrix_dh}
dh=\begin{pmatrix}
0&2&0\\
0&0&2
\end{pmatrix}\,.
\end{equation}
As a consequence of \eqref{eq:matrix_dh}, \eqref{eq:alpha_hopf} and \eqref{eq:a_wedge_da}, we obtain:

\begin{lemma}\label{lem:basic_prop_Hopf}
The Hopf map $h:\S^3\to \S^2$ is a smooth horizontally conformal map of constant dilation factor, namely
\begin{equation}\label{def:norm_of_grad_Hopf}
|dh|=2\sqrt{2}\quad \text{ and }\quad |h^\ast\omega_{\S^2}|=4\,,
\end{equation}	
and is of unit charge, i.e., $Q(h)=1$.
\end{lemma}

\smallskip

\section{The Hopf invariant for maps of finite Faddeev--Skyrme energy}\label{subsec:Hopf_inv}
\label{sec:degree}

In this section we prove some results concerning the definition and integrality of the Hopf invariant for maps with finite 
Faddeev--Skyrme  energy,  cf.\ \eqref{eq:energy_space}. The integrality 
 property  of the Hopf invariant  in this class was correctly  claimed in \cite[Section 3]{lin2004existence} but, as explained in \cite[\textsection 2.1.4]{auckly2010pontrjagin}, there is a gap in  the argument.  
In \cite[Corollary 5.2]{auckly2010pontrjagin} a 
 complete 
proof of integrality was given, but 
 the proof therein  is fairly long and not self-contained, as it relies on the earlier results by the authors in \cite{Auckly2005}. Here we present  a  new, more concise  
proof of this integrality which, although also not self-contained, we believe to be fairly transparent, and which contains ingredients which will  also  be important in the sequel.

In order to discuss how the notion of \textit{Hopf invariant} generalizes to maps with finite Faddeev--Skyrme energy, we first recall its definition and basic properties in the smooth setting, so let $u\in C^\infty(\S^3;\S^2)$. Since the volume-form $\omega_{\S^2}$ is closed and the exterior differential commutes with the pull-back, we have
\[d(u^*\omega_{\S^2})=u^*(d\omega_{\S^2})=0\,,\]
and in view of the fact that $H^2_{\rm{dR}}(\S^3)=0$ there exists $\beta\in \Omega^1(\S^3)$ such that $u^*\omega_{\S^2}=d\beta$. The \textit{Hopf invariant} or \textit{Hopf charge} of $u$, see \cite{whitehead1947expression}, is then given by the integral formula \eqref{eq:1_Hopf_inv}, and in particular $Q(u)\in \Z$.

\begin{remark}\label{rem:indep_hopf_inv}
The definition \eqref{eq:1_Hopf_inv} can easily be seen to be independent of the choice of the representative $\beta$. Indeed, if $\beta_1,\beta_2\in \Omega^1(\S^3)$ are such that $d\beta_1=d\beta_2=u^*\omega_{\S^2}$, using the fact that $H^1_{\rm{dR}}(\S^3)=0$ there exists $\psi\in C^\infty(\S^3)$ such that $\beta_1-\beta_2=d\psi$, and then the assertion follows from the identity
\[\beta_1\wedge d\beta_1=\beta_2\wedge d\beta_2+d(\psi d\beta_2)\,,\]
and Stokes' theorem. Moreover, the definition is also independent of the choice of the particular 2-form on $\S^2$, \textit{i.e.}, $\omega_{\S^2}$ can be replaced therein with any $\omega\in \Omega^2(\S^2)$ for which $\frac{1}{4\pi}\int_{\S^2}\omega=1$.  	
\end{remark}	

In order to generalize these results to maps with finite Faddeev--Skyrme energy, we will consider the function spaces
\begin{equation*}\label{eq:adm_S3_S3}
\mathcal A_3(\S^3):= \big\{v\in W^{1,2}(\S^3;\S^3)\colon\tp{cof}(d v)\in L^{3/2}(\S^3)\big\}\,,
\end{equation*} 
where $\tp{cof}(dv)$ is identified with the pull-back map $v^*\colon \Omega^{2}(\S^3)\to \Omega^{2}(\S^3)$. These function spaces  were introduced  
in the context of nonlinear elasticity theory in \cite{Ball1981a} and 
appear naturally  in our problem, since we have the following generalization of the results in \cite{Universite1992,Muller1994,Sverak1988}, due to \cite[Theorems 1.3 and 4.4]{Hamburger1999}:

\begin{proposition}\label{prop:commute}
If $v\in \mathcal A_3(\S^3)$, then $v^* d = d v^* \colon \Omega^{\ell}(\S^3)\to \Omega^{\ell+1}(\S^3)$ for every $\ell\in\{0,1,2\}$. Moreover, one can define \begin{equation}\label{eq:3-degree}
\deg(v):=\fint_{\S^3} v^*\omega_{\S^3}\,,
\end{equation}
where $\omega_{\S^3}$ is the volume-form on $\S^3$, and it holds that $\deg(v)\in \Z$.
\end{proposition}

We emphasize that all known proofs of Proposition \ref{prop:commute} yield the result \textit{directly} for maps in $\mc A_3(\mb S^3)$, \textit{i.e.}, the conclusion does not follow by approximation of maps in $\mc A_3(\mb S^3)$ with smooth maps, see already Remark \ref{rem:obst} below.

With these notations, and recalling the definition of $\mathcal U_{\mathcal{FS}}$ in \eqref{eq:energy_space},  the main result of this section is the following:	

\begin{theorem}\label{thm:hopfdeg}
Let $u\in \mathcal{U}_{\mathcal{FS}}$.
Then $u^*\omega_{\mb S^2}$ is 
an exact form, that is, there exists a 1-form $\beta\in W^{1,2}(\S^3;\bigwedge^1 T^*\S^3)$ so that
\begin{equation}\label{eq:ustarom_exact}
u^*\omega_{\S^2}=d\beta\,, \quad d^*\beta = 0\,. 
\end{equation}
Moreover, there exists a lift $\hat u\in \mathcal A_3(\mb S^3)$ such that $u= h\circ \hat u$ and 
\begin{equation}\label{eq:hopf_inv_degree}
Q(u) = \frac{1}{16\pi^2} \int_{\mb S^3} \beta \wedge d\beta = \deg(\hat u) \in \Z\,.
\end{equation}
\end{theorem}

In addition to Proposition \ref{prop:commute}, the following result due to  \cite{Bethuel1991} is the key ingredient in our proof of Theorem \ref{thm:hopfdeg}. 

\begin{theorem}\label{thm:bethuel}
A map $u\in W^{1,2}(\mb S^3;\mb S^2)$ is in the strong $W^{1,2}$-closure of $C^\infty(\mb S^3;\mb S^2)$ if and only if $u^*\omega_{\mb S^2}$ is distributionally closed.
\end{theorem}

\begin{proof} By \cite[Theorem~1]{Bethuel1991}, a map $u\colon \S^3\to\S^2$ is $W^{1,2}$-strongly approximable by smooth maps if and only if it pulls-back $2$-forms on $\S^2$ to (distributionally) closed 2-forms on $\S^3$.  Since any form in $\Omega^2(\mb S^2)$ can be represented as $f\omega_{\mb S^2}$ for some $f \in C^{\infty}(\S^2)$, it is then necessary and sufficient to check closedness of $u^*\omega_{\mb S^2}$.
\end{proof}

In order to apply Theorem \ref{thm:bethuel}, we will use the following:

\begin{lemma}\label{lem:closed}
If $u\in \mathcal{U}_{\mathcal{FS}}$, then $u^*\omega_{\mb S^2}$ is a distributionally closed 2-form.
\end{lemma}

\begin{proof}
This is detailed in \cite[Lemma 3.2]{auckly2010pontrjagin}. Fix $\lbrace j,k,\ell\rbrace\subset\lbrace 1,2,3\rbrace$.
Thanks to classical approximation by 
(localization and) mollification,
we see that $u^k,u^\ell\in W^{1,2}(\S^3;\R)$ satisfy
\begin{align*}
d(du^k\wedge du^\ell)=0\,,
\end{align*}
hence the  2-form 
$\alpha^{k\ell} \colon=du^k\wedge du^\ell$ 
is distributionally closed, and using again approximation by \PPP a \EEE mollification of $u^j\in W^{1,2}(\S^3;\R)$ and of the closed form $\alpha^{k\ell}\in L^2(\S^3;\bigwedge^2 T^*\S^3)$, we infer that
\begin{align*}
d(u^j\wedge\alpha^{k\ell})=d u^j\wedge \alpha^{k\ell}\,.
\end{align*}
Using this and \eqref{eq:omega_coordinates}, we obtain
\begin{equation}\label{eq:d_of_pullback}
d(u^*\omega_{\mb S^2}) =3 \, du^1 \wedge d u^2 \wedge d u^3\in L^1\Big(\S^3;\bigwedge\nolimits^3T^*\S^3\Big)\,,
\end{equation}
in the sense of distributions.

Since $u(x)\in \mb S^2$ for $\H^3$-a.e. $x\in \S^3$, \PPP which implies
\[u^1du^1+u^2du^2+u^3du^3=0\ \quad \text{for } \H^3\text{-a.e. } x\in \S^3\,,\] 
\EEE the vectors $(du^i(x))_{i=1}^3$ must be linearly dependent, hence the \PPP wedge product in the right hand side of \eqref{eq:d_of_pullback} \EEE vanishes.
\end{proof}

The next lemma, due to \cite{Hardt2003}, gives us a gauge for the lifts of $u$.

\begin{lemma}\label{lem:lift}	
Let $u\in \mathcal{U}_{\mathcal{FS}}$ and $\beta\in W^{1,2}(\S^3;\bigwedge^1 T^*\S^3)$ as in \eqref{eq:ustarom_exact}. Suppose that there is a lift $\hat u\in \mathcal{A}_{3}(\S^3)$ such that 
\begin{equation}\label{eq:hatu_alpha}
h\circ \hat u=u \ \text{ and }\ \beta = 2\hat u^*\theta\,,
\end{equation}
where, in the notation of \eqref{eq:alpha_exp} and \eqref{eq:alpha_hopf},
\begin{equation*}
h^*\omega_{\S^2}=2d\theta\,.
\end{equation*}
Then $\mathcal{H}^3$-a.e.\ on $\S^3$ it holds that
\begin{equation}\label{eq:lift}
|d \hat u|^2 =\frac 14   |\beta|^2 + \frac 14 |du|^2\,.
\end{equation}
Moreover, if $u\in C^\infty(\S^3;\S^2)$ then such a lift exists and is unique up to multiplication by a constant in $\mathbb S^1$ 
(that is, the action of 
$\S^1$ on $\S^3$ given in \textsection~\ref{subsec:horiz_conf}  by componentwise multiplication in $\C^2$).
\end{lemma}

\begin{proof}
The proof in the case of a smooth map $u$ can be found in \cite[Lemma 2.1]{Hardt2003}, 
and otherwise the conclusion of \eqref{eq:lift} can be repeated analogously, since the calculations therein can be carried over pointwise $\H^3$-a.e. on $\S^3$. The coefficients in \eqref{eq:lift} are actually different in \cite{Hardt2003} due to a different choice of normalization for the metric on $\mathbb S^2$ 
(which also changes the normalization constant in the integral formula \eqref{eq:1_Hopf_inv} for the Hopf degree).
For the readers' convenience we reproduce the proof of \eqref{eq:lift} with the metric used here. \EEE
In the orthonormal frame $\lbrace \tau_1,\tau_2,\tau_3\rbrace$ on $T\S^3$,
the $1$-form $\theta$ of \eqref{eq:alpha_exp} corresponds to the scalar product against $\tau_1$, cf. \eqref{eq:special_S3_frame}.
Thus, the condition $\beta =2\hat u^*\theta$  and \eqref{eq:local_param}, \eqref{eq:special_S3_frame}, imply 
\begin{equation}\label{eq:beta_norm}
\beta=2\tau_1(\hat u)\cdot d\hat u\implies|\beta|^2=4|\tau_1\cdot d\hat u|^2\,. 
\end{equation}	
Moreover, recalling the expression \eqref{eq:matrix_dh} of $dh$ in this orthonormal frame, the  chain rule gives
\begin{equation}\label{eq:du_chain_rule}
 du=dh\circ d\hat u=2(\tau_2\cdot d\hat u)f_1 + 2(\tau_3\cdot d\hat u)f_2\implies |du|^2=4|\tau_2\cdot d\hat u|^2 + 4 |\tau_3\cdot d\hat u|^2\,,
\end{equation}
 which, together with \eqref{eq:beta_norm} directly implies 
\[4|d\hat u|^2 =|\beta|^2+|du|^2\,,\]
as wished.
\end{proof}

\begin{proof}[Proof of Theorem \ref{thm:hopfdeg}]
Using Theorem \ref{thm:bethuel} and Lemma \ref{lem:closed} we can find a sequence of maps $(u_j)_{j\in \N}\subset C^\infty(\mb S^3;\mb S^2)$ such that 
\begin{equation}\label{eq:strong_W_1_2_approx}
u_j\rightarrow u \text{ strongly in } W^{1,2}(\S^3;\S^2) \text{ as }j\to\infty\,.
\end{equation} 
We can also fix forms 
$(\beta_j)_{j\in \N}\subset C^\infty(\S^3;\bigwedge^1 T^*\S^3)$
such that
\begin{align}\label{eq:dbeta_j}
d\beta_j= u^*_j\omega_{\mb S^2}\,, \qquad 
d^*\beta_j=0\,,
\qquad \fint_{\S^3}\beta_j=0\,.
\end{align} 
The last condition can be assumed without restriction after possibly subtracting a constant coefficient 1-form.
By standard elliptic estimates we have for instance the bound
\begin{equation}\label{eq:elliptic_estimate}
\|\beta_j\|_{L^{5/4}(\S^3)}\leq C \|u^*_j \omega_{\mb S^2}\|_{L^1(\S^3)}\,,
\end{equation} 
since $\frac 5 4 < \frac 3 2=: 1^*$. We note that the borderline $L^{3/2}$-estimate is actually also true, 
see \textit{e.g.} \cite{LanzaniStein2005} or \cite[Corollary~12]{BourgainBrezis2007}, at least in $\R^3$ instead of $\mb S^3$.
Invoking Lemma \ref{lem:lift},
we now fix lifts $(\hat u_j)_{j\in \N}\subset C^\infty(\S^3;\S^3)$
satisfying
\begin{equation}\label{eq:beta_j}
h\circ \hat u_j=u_j\,\ \text{and }\beta_j=2\hat u_j^*\theta\,, 
\end{equation}
cf.\ \eqref{eq:hatu_alpha}.
In view of \eqref{eq:am_gm-ineq} and \eqref{eq:strong_W_1_2_approx}-\eqref{eq:elliptic_estimate}, by passing to a non-relabelled subsequence, we obtain $\beta\in L^{5/4}\big(\S^3;\bigwedge\nolimits^1T^*\S^3\big)$
so that 
\begin{equation}\label{eq:weak_L_3_2_conv__forms}
\beta_j\rightharpoonup \beta \text{ weakly in } L^{5/4}\big(\S^3;\bigwedge\nolimits^1T^*\S^3\big) \text{ as }j\to\infty\,,
\end{equation}
and since convergence of averages is preserved under weak convergence, also
\begin{equation}\label{eq:average_beta_0}
\fint_{\S^3}\beta=0\,.
\end{equation}
By \eqref{eq:lift} applied to the triplet $(\hat u_j,\beta_j,u_j)$, together with \eqref{eq:strong_W_1_2_approx} and \eqref{eq:elliptic_estimate}, we have
\[\sup_{j\in\N}\|\hat u_j\|_{W^{1,5/4}(\S^3;\S^{ 3})}<+\infty\,,\]
and passing to a further non-relabeled subsequence, we obtain 
$\hat u\in W^{1,5/4}(\S^3;\S^{ 3})$
so that 
\begin{equation}\label{eq:weak_W_3_2_conv__lifts}
\hat u_j\rightharpoonup \hat u\text{ weakly in } W^{1,5/4}\big(\S^3;\S^{ 3}\big) \text{ as }j\to\infty\,.
\end{equation} 
Using now \eqref{eq:weak_L_3_2_conv__forms} and \eqref{eq:weak_W_3_2_conv__lifts}, we can pass to the limit in \eqref{eq:beta_j}, to infer that
\begin{equation}\label{eq:beta_equation}
h\circ \hat u=u\ \quad \H^3\text{-a.e.}\,,\qquad \beta = 2 \hat u^*\theta\,.
\end{equation} 
Indeed, in view of \eqref{eq:alpha_exp}, the forms $\hat u_j^*\theta$ can be written in coordinates as second-order polynomials of the form $\hat u_j\bullet d\hat u_j$, where the latter notation indicates that each monomial in the expression is of first order separately in $\hat u_j$ and $d\hat u_j$. In particular, since $\hat u_j\rightarrow \hat u$ strongly in $L^p(\S^3;\S^2)$ for every 
$ 1\leq p <+\infty
$ and $d\hat u_j\rightharpoonup d\hat u$ weakly in $L^{5/4}(\S^3;\bigwedge\nolimits^1(T^*\S^3\otimes \R^3))$, this particular product structure of weakly convergent objects and strongly convergent ones justifies \eqref{eq:beta_equation}. 

By passing to the limit in $L^1$ in \eqref{eq:dbeta_j}, using that \eqref{eq:omega_coordinates} and \eqref{eq:strong_W_1_2_approx} imply
\begin{equation*}
\|u_j^*\omega_{\S^2}-u^*\omega_{\S^2}\|_{L^1(\S^2)}\rightarrow 0 \ \text{ as } j\to\infty\,,
\end{equation*}
together with \eqref{eq:weak_L_3_2_conv__forms}, we also obtain
\begin{equation}\label{eq:dbeta_equation}
d\beta = u^*\omega_{\S^2}\,, \quad   d^*\beta=0\,.
\end{equation} 
Note however that $u\in \mathcal{U}_{\mathcal{FS}}$, which recalling \eqref{eq:energy_space} implies that $u^*\omega_{\S^2}\in L^{2}\big(\S^3;\bigwedge\nolimits^2T^*\S^3\big)$. Thus, \eqref{eq:dbeta_equation}, Gaffney's inequality \eqref{eq:gaffney} and \eqref{eq:average_beta_0} actually imply	 that $\beta\in W^{1,2}\big(\S^3;\bigwedge^1 T^*\S^3\big)$. The last assertion and \eqref{eq:lift} in turn imply that $\hat u \in W^{1,2}(\mb S^3;\mb S^3)$. 	
To conclude the proof it remains to verify that $\hat u \in \mathcal A_3(\mb S^3)$ and that \eqref{eq:hopf_inv_degree} holds true. 

As in the proof of \eqref{eq:lift}, 
 recall \eqref{eq:beta_norm}-\eqref{eq:du_chain_rule}, 
the identities $du=dh\circ d\hat u$ and $\beta=2\hat u^* \theta$
provide explicit expressions for the components of $d\hat u$ in terms of $du$ and $\beta$:
in the orthonormal frame $\lbrace \tau_1,\tau_2,\tau_3\rbrace$ on $T\S^3$ given by \eqref{eq:special_S3_frame}, 
the matrix of $d\hat u$ is given by
\begin{align*}
d\hat u
=\frac 12
\begin{bmatrix}
\beta(\tau_1) & \beta(\tau_2) & \beta(\tau_3)
\\
du(\tau_1)\cdot f_1
&
du(\tau_2)\cdot f_1
&
du(\tau_3)\cdot f_1
\\
du(\tau_1)\cdot f_2
&
du(\tau_2)\cdot f_2
&
du(\tau_3)\cdot f_2
\end{bmatrix}
=
\frac 12 
\begin{bmatrix}\beta\\du\cdot f_1 \\ du\cdot f_2 \end{bmatrix}
\,,
\end{align*}
where $\lbrace f_1,f_2\rbrace$ is the orthonormal frame on $T\S^2$
defined after \eqref{eq:dh_coord}.
Thus, using also \eqref{eq:modulus_ofu_omega}, we have
\begin{align*}
16|\tp{cof}(d\hat u)|^2
&
=4|d\hat u \wedge d \hat u|^2  
\\
&= 
|\beta \wedge (du\cdot f_1)|^2 +
|\beta\wedge (du\cdot f_2)|^2+ 
|(du\cdot f_1)\wedge (du\cdot f_2)|^2
\\
&
= 
|\beta \wedge du|^2 
+\frac 14 |du \wedge du|^2\,.
\end{align*}

Since $u\in \mathcal{U}_{\mathcal{FS}}$, hence $u^*\omega_{\mathbb S^2}\in L^2$, the last term is in $L^1$, cf.\ \eqref{eq:modulus_ofu_omega}. 
For the other term, we note that by \eqref{eq:dbeta_equation}, 
 \eqref{eq:average_beta_0}, \eqref{eq:gaffney} 
and the Sobolev embedding we obtain $\beta\in L^6\big(\S^3;\bigwedge\nolimits^1T^*\S^3\big)$, thus by Hölder's inequality 
$\beta\wedge du \in L^{3/2}(\S^3;\bigwedge\nolimits^2T^*\S^3)$, 
and we deduce that $\tp{cof}(d\hat u)\in L^{3/2}(\S^3)$.

Since $\hat u\in \mathcal{A}_{3}(\S^3)$,  we can apply Proposition \ref{prop:commute} to infer that $\deg(\hat u)$ is a well-defined integer and that $\hat u^*$ and $d$ commute. Hence, recalling \eqref{eq:alpha_hopf}, we have
\[u^*\omega_{\mb S^2} = \hat u^* h^*\omega_{\mb S^2} = \hat u^*(2d \theta) = d(2\hat u^*\theta)\,.\]
Thus, by \eqref{eq:a_wedge_da} and the fact that $\H^3(\S^3)=2\pi^2$, we obtain
\begin{align*}
 Q(u) &= \frac{1}{16 \pi^2}\int_{\mb S^3} (2\hat u^*\theta) \wedge d(2\hat u^*\theta)\\
 & =\frac{1}{4 \pi^2}\int_{\mb S^3} \hat u^*\theta \wedge \hat u^*(d \theta) =\frac{1}{4 \pi^2}\int_{\mb S^3}\hat u^*(\theta \wedge d \theta)= \frac{1}{2\pi^2}\int_{\mb S^3}\hat u^*\omega_{\S^3}=\deg(\hat u)\,,
\end{align*}
which concludes the proof.
\end{proof}

\begin{corollary}\label{cor:u_FS_comm_pullback}
If $u\in \mathcal{U}_{\mathcal{FS}}$, then $du^*=u^*d\colon \Omega^{1}(\S^2)\to \Omega^{2}(\S^3)$.
\end{corollary}

\begin{proof}
This property follows immediately from Theorem~\ref{thm:hopfdeg}, in particular the lifting identity $u=h\circ \hat u$, Proposition \ref{prop:commute}, and the smoothness of $h$, which also implies that $dh^*=h^*d$. 
\end{proof}

\begin{remark}\label{rem:obst}
There is a serious obstruction to having a more direct proof of Theorem \ref{thm:hopfdeg}, namely that it is not known whether one can approximate in the Faddeev--Skyrme energy a general map in  $\mathcal{U}_{\mathcal{FS}}$ by smooth maps. This is not known even if one ignores the image constraint on maps in $\mathcal {U}_{\mathcal{FS}}$, and we refer the reader to \cite{Hajasz2002} for some  partial results in the unconstrained case. We emphasize that from the above proof one cannot infer that the Hopf invariant of $u$ is obtained as the limit of the Hopf invariants of $(u_j)_{j\in \N}$ as in \eqref{eq:strong_W_1_2_approx}, since the latter sequence converges to $u$  only in 
too weak a sense for this purpose (and similarly for the liftings).
\end{remark}

We conclude this section with some easy but helpful consequences of Theorem \ref{thm:hopfdeg}. The first one is that  for maps of finite Faddeev--Skyrme energy the integral definition of the Hopf invariant is independent of the particular choice of the 2-form on $\S^2$, as long as its induced oriented area is fixed.

\begin{lemma}\label{lem:Hopf_indep_of_omega}
If $\omega_{1},\omega_{2}\in L^2(\S^2;\bigwedge^2T^*\S^2)$ are such that 
\begin{equation}\label{eq:omega_1_2_same_integral}
\int_{\S^2}\omega_1=\int_{\S^2}\omega_{2}
\end{equation}
and $v\in \mathcal{U}_{\mathcal{FS}}$ (cf. \eqref{eq:energy_space}) is such that $v^*\omega_{1},v^*\omega_{2}\in L^2(\S^3;\bigwedge^2T^*\S^3)$, then 
\begin{equation}\label{eq:Hopf_independ}
\int_{\S^3}\beta_1\wedge d\beta_1=\int_{\S^3}\beta_2\wedge d\beta_2\,,
\end{equation}
for any $\beta_1,\beta_2\in W^{1,2}(\mb S^3; \bigwedge^1 T^*\mb S^3)$
satisfying $d\beta_i=v^*\omega_{i}$ for $i=1,2$.
\end{lemma} 

\begin{proof}
Since the map $H^2_{\mathrm{dR}}(\S^2)\ni[\omega]\mapsto \int_{\S^2}\omega\in \R$ is a linear isomorphism, our assumption \eqref{eq:omega_1_2_same_integral} in particular implies that 
\[[\omega_{1}]=[\omega_{2}]\in H_{\mathrm{dR}}^2(\S^2)\,,\] 
thus there exists $\eta\in W^{1,2}(\S^2;\bigwedge^1T^*\S^2)$ such that 
\begin{equation*}
\omega_1=\omega_2+d\eta\,.
\end{equation*}
Using Corollary \ref{cor:u_FS_comm_pullback}, we can also calculate
\begin{align*}
d(\beta_1-\beta_2-v^*\eta)&=d\beta_1-d\beta_2-dv^*\eta=v^*\omega_1-v^*\omega_2-v^*d\eta\\
&=v^*\omega_1-v^*\omega_2-v^*(\omega_1-\omega_2)=0\,,
\end{align*}
and since $H^1_{\mathrm{dR}}(\S^3)=0$, for the 1-form $\beta_1-\beta_2-v^*\eta$ we can find $\gamma \in W^{1,2}(\S^3)$ such that
\begin{equation}\label{eq:fix_gauge}
\beta_1=(\beta_2+v^*\eta)+d\gamma\,.
\end{equation} 
In view of \eqref{eq:fix_gauge}, Remark \ref{rem:indep_hopf_inv} (which holds true even in this less regular setting) and \eqref{eq:forms_int_by_parts}, we have
\begin{align*}
\int_{\S^{3}}\beta_1\wedge d\beta_1&=\int_{\S^3}(\beta_2+v^*\eta)\wedge \big(d\beta_2+d(v^*\eta)\big)\\
&=\int_{\S^3}\beta_2\wedge d\beta_2+2\int_{\S^3}d\beta_2\wedge v^*\eta+\int_{\S^3}v^*\eta\wedge v^*d\eta\\
&=\int_{\S^3}\beta_2\wedge d\beta_2+2\int_{\S^3} v^*(\omega_2\wedge\eta)+\int_{\S^3}v^*(\eta\wedge d\eta)\\
&=\int_{\S^3}\beta_2\wedge d\beta_2\,,
\end{align*}
where in the last step we used the trivial fact that 
\[
\omega_2\wedge \eta, \eta\wedge d\eta\in L^1\big(\S^2;\bigwedge\nolimits^3T^*\S^2\big)=\{0\}\,,
\] 
hence the integrands in the last two summands vanish identically. This proves \eqref{eq:Hopf_independ}.
\end{proof}
As an immediate consequence of the previous lemma, we obtain the following \textit{degree formula}, for which we recall that for a map $\psi\in W^{1,2}(\S^2;\S^2)$ its \textit{degree} is defined, analogously to \eqref{eq:3-degree}, as the integral 
\begin{equation}\label{eq:2_degree}
\mathrm{deg}(\psi):=\fint_{\S^2}\psi^*\omega_{\S^2}\in \Z\,,	
\end{equation}
which in particular is continuous with respect to the strong-$W^{1,2}(\S^2;\S^2)$ topology.

\begin{lemma}\label{lem:degree_form}
Let $u\in \mathcal{U}_{\mathcal{FS}}$ and $\psi\in W^{1,2}(\S^2;\S^2)$ be such that $\psi^*\omega_{\S^2}\in L^2\big(\S^2;\bigwedge^2 T^*\S^2\big)$ and $\psi\circ u\in \mathcal{U}_{\mathcal{FS}}$. Then
\begin{equation}\label{eq:degree_formula}
Q(\psi\circ u)=\mathrm{deg}(\psi)^2\,Q(u)\,.
\end{equation}	
\end{lemma}

\begin{proof}	
By \eqref{eq:1_Hopf_inv} we have that 
\begin{equation}\label{eq:deg_comp}
16\pi^2\,Q(\psi\circ u)=\int_{\S^3}\beta_1\wedge d\beta_1\,,\ \text{where } d\beta_1=(\psi\circ u)^*\omega_{\S^2}=u^*(\psi^*\omega_{\S^2})=u^*\omega_1\,,
\end{equation}
where $\omega_1:=\psi^*\omega_{\S^2}\in L^2(\S^2;\bigwedge^2T^*\S^2)$. By \eqref{eq:2_degree} we observe that
\begin{equation*}
\int_{\S^2} \omega_1=\int_{\S^2}\psi^*\omega_{\S^2}= 4\pi\, \mathrm{deg}(\psi)\,. 
\end{equation*}	
Thus, setting $\omega_{2}:=\mathrm{deg}(\psi)\,\omega_{\S^2}\in \Omega^2(\S^2)$, we automatically have that 
\begin{equation*}\label{eq:omega_1_2_same_int}
\int_{\S^2}\omega_1=\int_{\S^2}\omega_2\,.
\end{equation*} 
Let now $\beta\in W^{1,2}\big(\S^3;\bigwedge^1T^*\S^3\big)$ be such that $d\beta=u^*\omega_{\S^2}$. Setting $\beta_2:=\mathrm{deg}(\psi)\,\beta$, we have
\begin{equation}\label{eq:beta_2}
d\beta_2=\mathrm{deg}(\psi)\,u^*\omega_{\S^2}=u^*(\mathrm{deg}(\psi)\,\omega_{\S^2})=u^*\omega_{2}\,.
\end{equation}
Thus, by \eqref{eq:deg_comp}-\eqref{eq:beta_2} and Lemma \ref{lem:Hopf_indep_of_omega}, we infer that 
\begin{equation*}
16\pi^2\,Q(\psi\circ u)=\int_{\S^3}\beta_1\wedge d\beta_1=	\int_{\S^3}\beta_2\wedge d\beta_2=\mathrm{deg}(\psi)^2\int_{\S^3}\beta\wedge d\beta=16\pi^2\mathrm{deg}(\psi)^2\,Q(u)\,,	
\end{equation*}
which shows \eqref{eq:degree_formula}. 
\end{proof}

\section{Relaxed form of the perturbed Faddeev--Skyrme energy}\label{sec:relaxed_formulation}
In order to prove Theorem \ref{thm:main}, and following the approach in \cite{isobe2008,riviere1998, riviere1998remark}, we consider a relaxation of $\mathcal{FS}_\rho$ to the space of closed differential forms. Let us first define, for a closed 2-form $\alpha$ on $\S^3$,
\begin{equation}\label{eq:Hopf_inv_forms}
Q(\alpha) := \frac{1}{16\pi^2} \int_{\mathbb S^3} \beta \wedge d\beta\,,\quad \alpha = d\beta\,.
\end{equation}
By Stokes' Theorem, this quantity, referred to as the \textit{Hopf invariant of the form $\alpha$},  is well-defined, \textit{i.e.},\ $Q$ depends on $\beta$ only through $d\beta$, see also the corresponding discussion in Remark \ref{rem:indep_hopf_inv} for the Hopf invariant of a map $u$. 
Introducing the space of admissible forms
\begin{equation*}\label{eq:adm_forms}
\mathcal{A}_{\mathcal{FS}}:=\left\{\alpha\in L^2\Big(\mathbb S^3; \bigwedge\nolimits^2 T^*\mathbb S^3\Big)\colon d\alpha=0\,, Q(\alpha)\PPP =1\EEE
\right\}\,,
\end{equation*}
we consider the relaxed energy $\mathcal{E}_\rho\colon \mathcal{A}_{\mathcal{FS}}\to [0,+\infty)$, defined by
\begin{equation}\label{eq:relaxed_FS_energy}
\mathcal{E}_\rho(\alpha):=\PPP 2\EEE 
\int_{\mathbb S^3}
|\alpha| +\PPP \rho^{-2}\EEE\int_{\mathbb S^3}|\alpha|^2\,.
\end{equation}

As an immediate consequence of the definitions of the energies and Lemma \ref{lem:dirichlet_energy_weak_conformality}, one obtains:
\begin{corollary}
\label{FS_maps_vs_forms}
 For  every $\rho>0$ and every $u\in \mathcal{U}_{\mathcal{FS}}$  it holds that
\begin{equation*}\label{eq:FS_maps_forms}
\mathcal{FS}_\rho(u) \geq \mathcal{E}_\rho(u^*\PPP\omega_{\S^2}\EEE)\,,
\end{equation*}
with equality if and only if $u$ is horizontally weakly conformal.
\end{corollary}

Recalling Proposition \ref{prop:spectrum}, we further introduce the following notation:
\begin{equation}\label{eq:defX}
	E_{0,1}^+:= \{\alpha \in E_0^+\colon Q(\alpha)=1\}\,.
\end{equation}
Using \eqref{eq:basis_for_E_0+}, it is straightforward to verify that 
\begin{equation}\label{eq:expression for E_01+}
E_{0,1}^+ =\big\lbrace c_1\omega^+_{0,1}+c_2\omega^+_{0,2}+c_3\omega^+_{0,3}\colon (c_1,c_2,c_3)\in \R^3 \ \text{with }c_1^2+c_2^2+c_3^2=16\big\rbrace\,.
\end{equation}

\BBB The next proposition for the relaxed energy $\E_\rho$ is the main ingredient in the proof of Theorem \ref{thm:main}. Compared with  
 \cite[Lemma~3.2]{isobe2008}, Proposition \ref{prop:stability} asserts global minimality of $E_{0,1}^+$, while the latter only gives minimality in a $C^0$-neighborhood of $E_{0,1}^+$.
\EEE

\begin{proposition}\label{prop:stability} 
\PPP For every $\rho\in(0,1]$ \BBB and \PPP$\alpha\in \mathcal{A}_{\mathcal{FS}}$, if 
\begin{equation}\label{eq:a_0_minim}
\alpha^+_{0,1}\in\underset{\tilde \alpha \in E_{0,1}^+}{\mathrm{argmin}}\,\|\alpha-\tilde\alpha\|_{L^2(\S^3)}\,,
\end{equation}
then 
\BBB
\begin{equation}\label{eq:_forms_minimality}
\mathcal{E}_\rho(\alpha)\geq \mathcal{E}_\rho(\alpha^+_{0,1})\,, 
\end{equation}
with equality if and only if $\alpha\in E_{0,1}^+$.
\end{proposition}
The rest of the section is devoted to the proof of Proposition \ref{prop:stability}. Before doing so, we present some auxiliary lemmata that will be essential for the proof. For $q=1,2$, let us introduce the auxiliary functionals $\I_q:\mathcal{A}_{\mathcal{FS}}\to [0,+\infty)$, defined as 
\begin{equation}\label{def:I_p}
\I_q(\alpha):=\int_{\mathbb S^3}|\alpha|^q\,.
\end{equation}
Observe that for $\alpha^+_{0,1}\in E^+_{0,1}$, setting 
\begin{equation}\label{eq:beta_0}
\beta_0:=\frac{1}{2}(*\alpha^+_{0,1}), \text{ we have } d\beta_0=\alpha^+_{0,1}\,.
\end{equation} 
The normalization condition for the Hopf invariant in \eqref{eq:defX}, together with the facts that $\alpha^+_{0,1}$ is a constant
coefficient 2-form, cf.\ \eqref{eq:expression for E_01+} and 
Remark~\ref{rem:structure_of_E_0}, that 
$$\beta_0\wedge d\beta_0=\tfrac{1}{2}(\ast\alpha^+_{0,1})\wedge \alpha^+_{0,1}=\tfrac{1}{2}|\alpha^+_{0,1}|^2\,,$$
and that $\H^3(\S^3)=2\pi^2$, imply that
\begin{equation}\label{eq:valueIp}
|\alpha^+_{0,1}|=4\,, \ \ \text{hence}\ \ \I_q(\alpha^+_{0,1})= 2\pi^2 4^q \ \ \text{for every } \alpha^+_{0,1}\in E^+_{0,1}\,.
\end{equation}

\PPP In \cite{speight2011} it is shown that elements in $E^+_{0,1}$ minimize the functional $\alpha \mapsto \I_2(\alpha)/Q(\alpha)$ on $\Omega^2(\mb S^3)$.
Actually, we can obtain an associated quantitative estimate, for both the eigenspaces $E_0^+$ and $E_0^-$, as follows.

\begin{proposition}\label{prop:faddeev}
The following estimates hold true:
\begin{itemize}
\item[$\rm{(i)}$] For every $\psi\in W^{1,2}(\S^3;\bigwedge^1 T^*\S^3)$ such that $d\psi\bot_{L^2}\,E_0^+\oplus E_0^-$, it holds that 
\begin{equation}\label{eq:psi_dpsi_est}
\Big|\int_{\S^3}\psi\wedge d\psi\Big|\leq \frac{1}{3}\int_{\S^3}|d\psi|^2\,.	
\end{equation}
\item[$\rm{(ii)}$] Let $\alpha\in L^2\big(\mb S^3;\bigwedge^2 T^*\mb S^3\big)$ be a  distributionally closed 2-form. Then 
\begin{equation*}\label{eq:quantitative_minim}
\I_2(\alpha)\mp32\pi^2 Q(\alpha)\geq \frac 13 \min_{\tilde \alpha^{\pm}_{0}\in E_0^\pm} \int_{\mb S^3} |\alpha-\tilde \alpha^\pm_{0}|^2 \,.
\end{equation*}
\end{itemize}
\end{proposition}

\begin{proof}
For (i), using Proposition \ref{prop:spectrum}, the $L^2$-orthogonality assumption, and since $H^2_{\rm{dR}}(\S^3)=0$, we can write 
\begin{equation*}
d\psi = \sum_{k=0}^\infty (\eta_k^++\eta_k^-)\,, \ \text{where } \eta_0^{\pm}=0\,,\  \eta^\pm_{k}=d\psi^\pm_{k}\in E_k^\pm\,,
\end{equation*}
where $\eta_k^\pm$ is the $L^2$-orthogonal projection of $d\psi$ onto $E_k^\pm$.
Since for every $\ell\in \N$ the 2-form $d\psi^{\pm}_\ell$ solves \eqref{eq:eigenval} with eigenvalue $\lambda^{\pm}_\ell:=\pm(\ell+2)$, using also Lemma \ref{lem:nulllag}(i),  we have that for $k,\ell\in \N$,
\begin{align}\label{eq:Hopf_orthogonality}
\begin{split}
\int_{\mb S^3} \psi^\pm_k\wedge d \psi^\pm_\ell &= \frac 1 {\lambda^{\pm}_\ell} \int_{\mb S^3} \psi_k^\pm \wedge d * d \psi_\ell^\pm =\frac 1 {\lambda^{\pm}_\ell} \int_{\mb S^3}(*d \psi_\ell^\pm) \wedge d \psi_k^\pm \\[2pt]
&= \frac 1 {\lambda^{\pm}_\ell} \int_{\mb S^3} \langle d \psi_\ell^\pm, d\psi_k^\pm\rangle=\frac{\delta^{k^\pm\ell^\pm}}{\lambda^\pm_k}\int_{\S^3}|d\psi_k^\pm|^2\,.
\end{split}
\end{align}
\BBB
In particular, \eqref{eq:Hopf_orthogonality} together with the fact that $|\lambda^\pm_{k}|\geq 3$ for $k\geq 1$, indeed implies 
\begin{align*}
\Big|\int_{\S^3}\psi\wedge d\psi\Big|&=	\Big|\sum_{k,\ell\geq 1}\int_{\S^3}\psi_k^\pm\wedge d\psi_\ell^\pm\Big|\leq	\sum_{k\geq 1}\frac{1}{|\lambda^\pm_k|}\int_{\S^3}|d\psi_k^\pm|^2\leq \frac{1}{3}\int_{\S^3}|d\psi|^2\,,
\end{align*}
with equality holding for $\psi_1^\pm$.
\PPP

For the proof of (ii), consider a closed $2$-form $\alpha$, thus $\alpha=d\beta$ for some $\beta\in W^{1,2}(\S^3;\bigwedge^1 T^*\S^3)$ with $d^*\beta=0$.
The lowest eigenvalue of the Hodge Laplacian $\Delta:=dd^*+d^\ast d$ on closed $2$-forms is $4$, and the corresponding eigenspace is given by  
\begin{align*}
\ker d \cap \ker(\Delta-4)= E_0^+ \oplus  E_0^-\,,
\end{align*}
cf.\ \cite[Section 3]{isobe2008} and \eqref{eq:eigenspaces}. 
Therefore, we can ($L^2$-orthogonally) decompose
\begin{equation}\label{eq:beta_decomp}
*\beta= \alpha_0^+ + \alpha_0^- + d\varphi\,, \ \text{where } \ d*\alpha_0^\pm=\pm 2 \alpha_0^\pm\,, \ \ d\varphi \perp_{L^2} E_0^+ \oplus  E_0^-\,,
\end{equation}
which, by taking the Hodge-dual and then the exterior differential results in 
\begin{equation}\label{eq:alpha_decomp}
\alpha=2\alpha_0^+ -2\alpha_0^- + d\psi, \ \text{where } \ \psi := *d\varphi\,, \ \ d\psi \perp_{L^2} E_0^+ \oplus  E_0^-\,.
\end{equation}
Using \eqref{eq:beta_decomp} and \eqref{eq:alpha_decomp}, we find
\begin{equation}\label{eq:I_2_exp}
\I_2(\alpha)=4\int_{\mathbb S^3} |\alpha_0^+|^2+ 
4\int_{\mathbb S^3}|\alpha_0^-|^2+\int_{\mathbb S^3} |d\psi|^2\,,
\end{equation}
and 
\begin{equation}\label{eq:Q_expansion}
Q(\alpha)=\frac{1}{16\pi^2}\int_{\mathbb S^3} \langle *\beta, \alpha\rangle\\
=\frac{1}{8\pi^2}\int_{\mathbb S^3}|\alpha_0^+|^2 -\frac{1}{8\pi^2}\int_{\mathbb S^3} |\alpha_0^-|^2+\frac{1}{16\pi^2}\int_{\mathbb S^3}\psi\wedge d\psi\,.
\end{equation}
We now proceed to give a lower bound on the last term in the right hand side of \eqref{eq:Q_expansion}. 
 Combining \eqref{eq:I_2_exp} and \eqref{eq:Q_expansion}, using the estimate \eqref{eq:psi_dpsi_est} and again \eqref{eq:alpha_decomp}, we deduce
\begin{align*}
\begin{split}	
\I_2(\alpha)\mp32\pi^2 Q(\alpha)&=8\int_{\mathbb S^3}|\alpha_0^\mp|^2 + \int_{\mathbb S^3}|d\psi|^2 \mp 2\int_{\mathbb S^3}\psi\wedge d\psi\\
&\geq 8\int_{\mathbb S^3}|\alpha_0^\mp|^2 +(1-\tfrac{2}{3}) \int_{\mathbb S^3}|d\psi|^2\\
&\geq \frac{1}{3}\int_{\S^3}|\alpha\mp2\alpha_0^{\pm}|^2 \geq \frac 13 \min_{\tilde \alpha^\pm_{0}\in E_0^\pm}\int_{\mathbb S^3}|\alpha - \tilde \alpha^\pm_{0}|^2\,,
\end{split}
\end{align*}
concluding the proof.
\end{proof}\EEE

The next lemma gives us a precise expansion for the energy $\E_\rho$ around elements in $E^+_{0,1}$. 

\begin{lemma}\label{lem:expE}
Given $\varphi \in W^{1,2}(\mb S^3; \bigwedge^1 T^*\mb S^3)$, let $\alpha:=\alpha^+_{0,1}+d\varphi$, where $\alpha^+_{0,1}\in E^+_{0,1}$. Then, 
\begin{align}\label{eq:expansionE}
\E_\rho(\alpha)-\E_\rho(\alpha^+_{0,1})&=\PPP \Big(\frac{1}{2}+\frac{2}{\rho^2}\Big)\int_{\S^3}\langle \alpha^+_{0,1},d\varphi \rangle+\frac{1}{\rho^2}\int_{\S^3}|d\varphi|^2+\mathcal{G}_1(\alpha^+_{0,1},d\varphi)\EEE\,,
\end{align}
where we have set
\begin{equation}\label{eq:G_1_formula}
\mathcal{G}_1(\alpha^+_{0,1},d\varphi):=\int_{\mathbb S^3} \bigg\lbrace
|\alpha^+_{0,1}+d\varphi| -|\alpha^+_{0,1}| - \Big\langle\frac{\alpha^+_{0,1}}{|\alpha^+_{0,1}|}, d\varphi\Big\rangle
\bigg\rbrace
\geq 0\,.
\end{equation} 
\end{lemma}

\begin{proof}
Since $H^2_{\rm{dR}}(\S^3)=0$, by the $L^2$-Hodge decomposition we can assume that $\varphi$ satisfies
\begin{equation*}
d^* \varphi=0\,.
\end{equation*}
Using the notation in \eqref{def:I_p},  the relaxed energy in \eqref{eq:relaxed_FS_energy} can be written as
\begin{equation}\label{eq:rewriting_the_energy}
\mathcal{E}_{\rho}(\alpha)=\PPP 2
\I_1(\alpha)+\rho^{-2}\I_2(\alpha)\EEE\,.	
\end{equation} 
We first expand the two quantities appearing in \eqref{eq:rewriting_the_energy} separately, obtaining:
\begin{align}\label{eq:three_terms_expansion}
\begin{split}	
{\rm(i)}&\quad \I_1(\alpha)	=\I_1(\alpha^+_{0,1})+\I_1'(\alpha^+_{0,1})[d\varphi] +\mathcal{G}_1(\alpha^+_{0,1},d\varphi)\,,\\[2pt]
{\rm(ii)}&\quad \I_2(\alpha)=\I_2(\alpha^+_{0,1})+\I_2'(\alpha^+_{0,1})[d\varphi] + \int_{\mathbb S^3} |d\varphi|^2\,,
\end{split}
\end{align}
where we have used the standard notation 
$$\I'_{q}(\alpha^+_{0,1})[d\varphi]:=\frac{d}{dt}\Big|_{t=0}\I_{q}(\alpha^+_{0,1}+td\varphi)\,,$$
and similarly for $Q'(\alpha^+_{0,1})[d\varphi]$, \PPP which we will also compute here for later purposes\EEE. 
Note that \eqref{eq:three_terms_expansion}(i) follows simply by the definition of $\mathcal{G}_1$ in \eqref{eq:G_1_formula}, and \eqref{eq:three_terms_expansion}(ii) by a simple expansion of the quadratic energy $\I_2$.

To proceed further we note that, by the definition of $\beta_0$ in  \eqref{eq:beta_0} and Lemma \ref{lem:nulllag}(i), we have 
\begin{equation}\label{eq:Q'}
16 \pi^2 Q'(\alpha^+_{0,1})[d\varphi]=\int_{\mathbb S^3}\beta_0\wedge d\varphi + \int_{\mathbb S^3}\varphi\wedge d\beta_0
=2\int_{\mathbb S^3}\beta_0\wedge d\varphi
=\int_{\mathbb S^3} \langle \alpha^+_{0,1}, d\varphi\rangle\,.
\end{equation}
 Then, a simple calculation using \eqref{eq:valueIp} and \eqref{eq:Q'} (with a general test form in the place of $d\varphi$) reveals that, for $q=1,2$,
\begin{equation}\label{eq:criticality}
\I_q'(\alpha^+_{0,1}) = \frac{q}{2}\I_q(\alpha^+_{0,1})Q'(\alpha^+_{0,1})\,;
\end{equation}
this last equation, together with $Q(\alpha^+_{0,1})=1$, expresses the fact that $\alpha^+_{0,1}$ is  a critical point of  
the functionals $\alpha\mapsto Q^{-q/2}(\alpha) \I_q(\alpha)$. We now analyze each of these two functionals separately: using  \eqref{eq:valueIp} and \eqref{eq:three_terms_expansion}-\eqref{eq:criticality}, we find that
\begin{align}\label{eq:two_expansions}
\begin{split}		
\rm{(i)}\quad & \PPP
\I_1(\alpha)-\I_1(\alpha^+_{0,1})
=\PPP \frac{1}{4}\int_{\S^3}\langle \alpha^+_{0,1},d\varphi \rangle+\mathcal{G}_1(\alpha^+_{0,1},d\varphi)\,,
\\[5pt]
\rm{(ii)}\quad &\PPP \I_2(\alpha)
- \I_2(\alpha^+_{0,1})
=\PPP \int_{\S^3}|d\varphi|^2+ 2\int_{\S^3}\langle \alpha^+_{0,1},d\varphi \rangle\EEE\,.
\end{split}		
\end{align}
Collecting \eqref{eq:relaxed_FS_energy}, and the expansions \eqref{eq:two_expansions}(i)-(ii), we find precisely \eqref{eq:expansionE}. 
\end{proof}

\begin{proof}[Proof of Proposition \ref{prop:stability}]
\PPP In view of \EEE \eqref{eq:basis_for_E_0+} 
 and  \eqref{eq:expression for E_01+}, we can identify $E_0^+\cong \R^3$ and $E^+_{0,1}\cong \mb S_4^2$.  Using \eqref{eq:L_2_orthogonal_Hodge} it is immediate that \PPP\eqref{eq:a_0_minim} \EEE also implies 
\begin{equation*}
\PPP\alpha^+_{0,1} \in \underset{\tilde \alpha \in E_{0,1}^+}{\mathrm{argmin}}\EEE\,\big\|\PPP \alpha_0^+-\tilde\alpha\big\|_{L^2(\S^3)}\,,
\end{equation*}	
where \PPP $\alpha_0^+:=P_{E_0^+}\alpha$ is the $L^2$-projection of $\alpha$ into $E_0^+$.
\EEE
Thus, $(\alpha_0^+-\alpha^+_{0,1})$ is in the normal space to $E^+_{0,1}\cong \mb S_4^2\subset \R^3$ at $\alpha^+_{0,1}$, \textit{i.e.},
\begin{equation*}\label{eq:P_E_0+alpha_paral}
\PPP\alpha_0^+\EEE=(1+t)\alpha^+_{0,1}\, \ \text{for some } t\in \R\,.
\end{equation*}  
Actually, since $d(\alpha-\alpha^+_{0,1})=0$ and $H^2_{\rm{dR}}(\S^3)=0$, there exists $\varphi \in W^{1,2}\big(\mb S^3; \bigwedge^1 T^*\mb S^3\big)$ with 
\begin{equation}\label{eq:dphi} 
\alpha=\alpha^+_{0,1}+d\varphi\ \text{ and }\ d^\ast\varphi=0\,.
\end{equation} 
Recalling \eqref{eq:beta_0}, we can write $\alpha = \alpha^+_{0,1} + d \varphi = (1+t)\alpha^+_{0,1} + (d \varphi-t \,d \beta_0)$, or equivalently, 
\begin{equation}\label{eq:psi}
\PPP \alpha= \alpha_0^+ +d \psi\,, \quad \text{for } \psi:=\varphi-t \beta_0\,,\EEE
\end{equation}
where we note that \PPP $d \psi \in (E_0^+)^\bot$\EEE, the latter intended in the $L^2$-sense:  hence, by \eqref{eq:psi} and using also \eqref{eq:valueIp}, we have \EEE
\begin{equation}\label{eq:L2orthog}
\int_{\S^3}\langle \alpha^+_{0,1}, d \psi\rangle=0\PPP \iff t=\frac{1}{32\pi^2}\int_{\S^3}\langle \alpha^+_{0,1},d\varphi\rangle\,.
\end{equation}
\PPP Moreover,  since $Q(\alpha)=Q(\alpha^+_{0,1})=1$, 
by \eqref{eq:forms_int_by_parts} we have
\begin{align*}
&1=Q(\alpha^+_{0,1}+d\varphi)=Q(\alpha^+_{0,1})+Q(d\varphi)+\frac{1}{8\pi^2}\int_{\S^3}\beta_0\wedge d\varphi=1+Q(d\varphi)+2t\,,
\end{align*}
where in the last equality we used again \eqref{eq:beta_0}, \eqref{eq:valueIp}, and  \eqref{eq:L2orthog}. Thus,
\begin{equation}\label{eq:t_equation}
\frac{1}{32\pi^2}\int_{\S^3}\langle \alpha^+_{0,1},d\varphi\rangle=t=-\frac 1 2 Q(d\varphi)\,.
\end{equation}
Note further that by \eqref{eq:a_0_minim}, and since $-\alpha^+_{0,1}\in E_{0,1}^+$, we have
\begin{equation*}
\|\alpha-\alpha^+_{0,1}\|_{L^2(\S^3)}\leq \|\alpha+\alpha^+_{0,1}\|_{L^2(\S^3)}\iff \int_{\S^3}\langle\alpha,\alpha^+_{0,1}\rangle\geq 0\,,
\end{equation*}
which, by \eqref{eq:dphi} is equivalent to the simple bound
\begin{equation}\label{eq:t_Q_simple_bound}
t\geq -1\iff Q(d\varphi)\leq 2	\,.
\end{equation}
The definition of $\psi$ in \eqref{eq:psi}, together with \eqref{eq:L2orthog}-\eqref{eq:t_Q_simple_bound} further yield the relation
\begin{align}\label{eqQ_d_phi_Q_psi}
\begin{split}	
&Q(d\varphi)=Q(d\psi)+t^2Q(\alpha^+_{0,1})=Q(d\psi)+\frac{1}{4}[Q(d\varphi)]^2\\[3pt]
\implies& Q(d\varphi)=2(1-\sqrt{1-Q(d\psi)})\,,
\end{split}
\end{align}
since \eqref{eq:t_Q_simple_bound} and the first line in the above identity imply that $Q(d\psi)\leq 1$.

\PPP
We can now apply Lemma \ref{lem:expE} and use \eqref{eq:G_1_formula}, \eqref{eq:psi}, and \eqref{eq:t_equation}, to find 
\begin{align}\label{eq:expE1}
\begin{split}
\E_\rho(\alpha)-\E_\rho(\alpha^+_{0,1})
&
\geq \rho^{-2}\int_{\mathbb S^3}|d\varphi|^2-16\pi^2\Big(\frac{1}{2}+\frac{2}{\rho^2}\Big)Q(d\varphi)\\
&= 
\rho^{-2}\Big(\mathcal{I}_2(d\varphi)-32\pi^2Q(d\varphi)\Big)-8\pi^2Q(d\varphi)\\
&
=
\rho^{-2}\Big(\mathcal{I}_2(d\psi)-32\pi^2Q(d\psi)\Big)
-16\pi^2\big(1-\sqrt{1-Q(d\psi)}\big)\,.
\end{split}
\end{align}
\BBB The last equality follows from \eqref{eqQ_d_phi_Q_psi}, the fact that
\[\mathcal I_2(d\varphi)=\mathcal I_2(d\psi)+t^2\mathcal I_2(\alpha_{0,1}^+)\] 
by the orthogonality in \eqref{eq:L2orthog},  and that 
\[\mathcal I_2(\alpha_{0,1}^+)=32\pi^2Q(\alpha^+_{0,1})=32\pi^2\,,\]
cf.\ \eqref{eq:valueIp}. \PPP
We rewrite \eqref{eq:expE1} as
\begin{align*}
\E_\rho(\alpha)-\E_\rho(\alpha^+_{0,1})
&
\geq 
\Big(\frac{1}{\rho^2}-\frac 12\Big)\Big(\mathcal{I}_2(d\psi)-32\pi^2Q(d\psi)\Big)
-16\pi^2\Big(1-\sqrt{1-Q(d\psi)}-\frac 12 Q(d\psi)\Big)
\notag\\
&
\quad
+\frac 12 
\Big(\mathcal{I}_2(d\psi)-48\pi^2Q(d\psi)\Big)\,.
\end{align*}
Note that in view of 
\eqref{eq:Hopf_inv_forms}-\eqref{eq:psi_dpsi_est} 
combined with decomposing $d\psi\in (E_0^+)^\perp$ according to the orthogonal sum $(E_0^+)^\perp
=(E_0^+\oplus E_0^-)^\perp\oplus E_0^{-}$,
the last line in the above inequality is nonnegative, hence 
we deduce
\begin{align}
\E_\rho(\alpha)-\E_\rho(\alpha^+_{0,1})
&
\geq 
\frac{2\rho^{-2}-1}{2}\Big(\mathcal{I}_2(d\psi)-32\pi^2Q(d\psi)\Big)
\nonumber
\\
&
\quad
-16\pi^2\Big(1-\sqrt{1-Q(d\psi)}-\frac{Q(d\psi)}{2}\Big)\,,
\label{eq:deficit_bound_psi}
\end{align}
which can be equivalently written as, \begin{align}\label{eq:deficit_almost_final}
\begin{split}
\frac{\E_\rho(\alpha)-\E_\rho(\alpha^+_{0,1})}{16\pi^2}
&
\geq 
(2\rho^{-2}-1)\Big(\frac{\mathcal{I}_2(d\psi)}{32\pi^2}-Q(d\psi)\Big)-\frac{1}{2}|Q(d\psi)|\\
&
\quad+\frac{1}{2}|Q(d\psi)|-\Big(1-\sqrt{1-Q(d\psi)}-\frac{Q(d\psi)}{2}\Big)
\,.
\end{split}
\end{align}
By the elementary fact that the smallest value of $\lambda>0$ such that the function
\begin{align*}
(-\infty,1]\ni s\mapsto \lambda |s|-\Big(1-\sqrt{1-s}-\frac s2\Big)\,,
\end{align*}
is nonnegative is $\lambda=1/2$, \eqref{eq:deficit_almost_final} together with \eqref{eq:Hopf_inv_forms} and \eqref{eq:psi_dpsi_est} again, finally yield, 
\begin{equation*}
\PPP\E_\rho(\alpha)-\E_\rho(\alpha^+_{0,1})\BBB\geq 
(2\rho^{-2}-1)\Big(\frac{\mathcal{I}_2(d\psi)}{\PPP 2\BBB}-\frac{\mathcal{I}_2(d\psi)}{\PPP 3\BBB}\Big)-\frac{\mathcal{I}_2(d\psi)}{\PPP 6\BBB}=
\Big(\frac{\rho^{-2}-1}{\PPP 3\BBB}\Big)
\mathcal{I}_2(d\psi) \geq 0
\,,
\end{equation*}
the latter inequality holding for $0<\rho \leq 1$. This proves the desired minimality \PPP \eqref{eq:_forms_minimality} \BBB of $\alpha_{0,1}^+$.

For the characterization of the equality cases, the equality $\E_\rho(\alpha)=\E_\rho(\alpha^+_{0,1})$ implies that equality holds throughout the above chain of inequalities. In particular, it holds in the first inequality in \eqref{eq:expE1}, which implies that $\mathcal G_1(\alpha_{0,1}^+,d\varphi)=0$. Thus the integrand in  \eqref{eq:G_1_formula} vanishes a.e.\ and so, by strict convexity of the Euclidean norm, we find a measurable $\lambda\colon \mb S^3\to \R$ such that
$d\varphi(x) = \lambda(x)\alpha^+_{0,1}\,.$
Taking a weak exterior derivative in this equation, since $\alpha^+_{0,1}$ is closed, we see that $\lambda$ is constant, say $\lambda(x)=t$ for all $x\in \mb S^3$, thus by \eqref{eq:dphi} we have $\alpha=(1+t) \alpha_{0,1}^+\in E_0^+$. But by assumption $Q(\alpha)=1$, so $\alpha\in E_{0,1}^+$, as wished; in particular, we must have \PPP $(1+t)^2=1$, hence $t=0$, 
since the case $t=-2$ is excluded by \eqref{eq:t_Q_simple_bound}.
\end{proof}
\EEE 
\begin{remark}[\textit{On the threshold for linear stability}]\label{rem:linear_stability}
A similar spectral analysis as the one in the proof of Proposition \ref{prop:faddeev} had already been used by Isobe in \cite{isobe2008} to prove that 
\[\textit{ the Hopf map } h \textit{ is a stable critical point for } \mathcal{FS}_\rho \textit{ for every }\rho\in(0,\sqrt{2})\,,\]
and by continuity weakly stable also for $\rho_\ast=\sqrt{2}$, cf.\ the first assertion of Theorem 1.1, Proposition 4.1 and Corollary 4.1 therein. The role of the threshold $\rho_*=\sqrt{2}$ for the linear stability of $h$ is also revealed by our calculations in the proof Proposition \ref{prop:stability}, for instance in \eqref{eq:deficit_bound_psi}. Indeed, by Proposition \ref{prop:faddeev} and that
\[1+\theta\geq 2\sqrt \theta\ \quad \forall \theta\geq 0\,,\]
the latter applied for $\theta:=1-Q(d\psi)$, if $\rho>\sqrt{2}$ the right hand side of \eqref{eq:deficit_bound_psi} becomes negative, leading to an insufficient lower bound for the conclusion.
\end{remark}

As an immediate consequence of Proposition \ref{prop:stability}, we obtain the minimality of the Hopf map for the $\mathcal{FS}_\rho$ in the claimed range.

\begin{corollary}\label{corol:minimality}
For every $\rho\in (0, 1]$ and $u\in \mathcal{U}_{\mathcal{FS}}$ with $Q(u)=1$, it holds that
\begin{equation*}\label{eq:claimed_minimality_hopf}
\mathcal{FS}_\rho(u) \geq \mathcal{FS}_\rho(h)\,,
\end{equation*}	
\BBB
with equality if and only if $u^*\omega_{\mb S^2}\in E_{0,1}^+$.
\EEE
\end{corollary}
\begin{proof}
\BBB Indeed for every $\rho\in (0,1\PPP ]\BBB$ and $u$ as in the statement, by Corollary \ref{FS_maps_vs_forms}, \eqref{eq:alpha_hopf}, and Proposition \ref{prop:stability},  we obtain
\begin{equation}\label{eq:minimality}
\mathcal{FS}_\rho(u) -\mathcal{FS}_\rho(h) \geq \E_\rho(u^*\omega_{\S^2})-\E_\rho(h^*\omega_{\S^2})
\geq 0\,.
\end{equation}
The characterization of the equality cases also follows from Proposition \ref{prop:stability}.
\end{proof}
\EEE

\section{Proof of uniqueness in Theorem \ref{thm:main}}\label{sec:proof_main}
The uniqueness statement in Theorem \ref{thm:main} will follow from the following:

\begin{proposition}\label{prop:weakly_conf_rigid}
Let $u\in  \mc U_{\mathcal FS}$ be a horizontally weakly conformal map such that $Q(u)=1$ and $u^*\omega_{\mb S^2} = h^*\omega_{\mb S^2}$. Then $u=h\circ R$ for some $R\in SO(4)$.
\end{proposition}

\begin{proof}
We first claim that 
\begin{equation}\label{eq:du_tau_1_0}
du(\tau_1)=0 \quad \H^3\text{-a.e. on } \S^3\,,	
\end{equation}
where $\tau_1$ is the fundamental vertical vector field for the Hopf fibration defined in \eqref{eq:special_S3_frame}.
By the assumptions, $\H^3$-a.e. on $\S^3$ the frame $\{\tau_1,\tau_2,\tau_3\}$ satisfies for $i,j\in\{1,2,3\}$
\begin{equation}\label{eq:du_orthog}
\omega_{\S^2}(du(\tau_i),du(\tau_j))=u^*\omega_{\S^2}(\tau_i,\tau_j)=h^*\omega_{\S^2}(\tau_i,\tau_j)=\omega_{\S^2}(dh(\tau_i),dh(\tau_j))\,.
\end{equation} 
Since $\{dh(\tau_2),dh(\tau_3)\}$ span $T\S^2$, cf.\ \eqref{eq:dh_coord}, \eqref{eq:du_orthog} implies that 
\begin{equation}\label{eq:dh_span}
\{du(\tau_2),du(\tau_3)\} \ \text{are linearly independent and span } T\S^2\,.
\end{equation}
Using again \eqref{eq:dh_coord} and \eqref{eq:du_orthog}, this time for $j=2,3$, we have
\begin{equation*}
\omega_{\S^2}(du(\tau_1),du(\tau_j))=\omega_{\S^2}(dh(\tau_1),dh(\tau_j))=0\,,
\end{equation*} 
which by \eqref{eq:dh_span} and the non-degeneracy of $\omega_{\S^2}$ implies \eqref{eq:du_tau_1_0}. 

We can now define $\psi\in W^{1,2}(\S^2;\S^2)$ with the property that 
\begin{equation*}\label{eq:psi_prelift}
u=\psi\circ h\,.
\end{equation*}
Indeed, since $h\colon \S^3\to \S^2$ is onto, for every $y\in \mb S^2$ we can find $x\in \mb S^3$ such that $h(x)=y$ and then we set $\psi(y) := u(x)$. This is well-defined, \textit{i.e.}, for $\H^3$-a.e.\ $x\in \S^3$, if $x'\in \S^3$ is such that $h(x)=h(x')$, so that $x,x'$ belong to the same $\S^1$-fiber of $h$, then $u(x)=u(x')$.  For this, observe that by the general slicing properties of Sobolev functions, $u$ is  absolutely continuous  on a generic  $\S^1$-fiber, and the claim follows from \eqref{eq:du_tau_1_0}. That $\psi\in W^{1,2}(\S^2;\S^2)$ follows trivially by the fact that $u\in W^{1,2}(\S^3;\S^2)$ and that $h\colon \S^3\to \S^2$ is a $C^\infty$-submersion, cf.\ \cite[Theorem 4.29]{lee2003smooth} for the smooth case. In particular, we can pick smooth local coordinates $(\tilde x_1, \tilde x_2, \tilde x_3)$ on $\mb S^3$ such that $h(\tilde x) = (\tilde x_1, \tilde x_2)\in \mb S^2$ and we get the local expression $\psi(\tilde x_1, \tilde x_2) = u(\tilde x_1, \tilde x_2,\tilde x_3)$,  hence $\psi$ has the same regularity as $u$.

Since $h$ is smooth,  by the chain rule for the weak derivatives we have 
\[du(x)= d\psi(h(x)) \circ dh(x) \text{ for } \H^3\text{-a.e.\ } x\in \S^3\,,\] 
which together with the weak horizontal conformality \eqref{eq:horiz_weak_conf_alternative} of $u$ and $h$  gives
\begin{align*}\label{eq:conf_psi}	
d\psi(h(x))\circ dh(x)\circ(dh(x))^t\circ(d\psi(h(x)))^t=\frac{1}{4}|d\psi(h(x))|^2|dh(x)|^2\,\mathrm{Id}_{T_{\psi(h(x))}\S^2}\nonumber
\end{align*}	
or, equivalently, 
$$d\psi(h(x))\circ(d\psi(h(x)))^t=\frac{1}{2}|d\psi(h(x))|^2\,	\mathrm{Id}_{T_{\psi(h(x))}\S^2}\,,$$
\textit{i.e.}, $\psi\in W^{1,2}(\S^2;\S^2)$ is a weakly conformal map.  Moreover, our assumptions, \eqref{def:norm_of_grad_Hopf} and the weak conformality of $\psi$ directly imply that
\[4=|h^*\omega_{\mb S^2}| = |u^* \omega_{\mb S^2}| = |h^* \psi^* \omega_{\mb S^2}| = \frac{|d\psi|^2}{2} |h^*\omega_{\mb S^2}|=2|d\psi|^2\,,\]
 \textit{i.e.}, $d\psi$ is  $\mathcal{H}^2$-a.e.\ an isometry with $|d\psi|^2 = 2$. 
Using Lemma \ref{lem:degree_form}, which is applicable since $|\psi^*\omega_{\S^2}|=\frac{|d\psi|^2}{2}=1$ and $\psi\circ h=u\in \mathcal{U}_{\mathcal{FS}}$, we also obtain 
\[1 = Q(u) = Q(\psi \circ h) = \deg(\psi)^2 Q(h) = \deg(\psi)^2\,,\]
\textit{i.e.}, we conclude that 
\begin{equation*}\label{eq:conclusion_for_psi}
\psi\in W^{1,2}(\S^2;\S^2) \text{ is such that } d\psi^td\psi=\mathrm{Id}_{T\S^2}\, \text { with } \mathrm{deg}(\psi)=\pm 1\,.
\end{equation*}
It follows that $\psi\equiv O\mathrm{id}_{\S^2}$ for some $O\in O(3)$, as can easily be inferred from Liouville's theorem, see for instance \cite[Theorem 1.1]{hirsch2022note}. Finally,  using the equivariance of the Hopf map  exactly as at the end of the proof of \cite[Proposition 4.2]{isobe2008}, there exists $R\in \tp{O}(4)$ such that 
\begin{equation*}\label{eq:u_conclusion} 
u=\psi \circ h = h \circ R\,,
\end{equation*}
and since $Q(u)=Q(h)=1$, by applying Theorem \ref{thm:hopfdeg}  we actually deduce that $R\in \tp{SO}(4)$, as desired. 
\end{proof}

\begin{proof}[Completion of proof of Theorem \ref{thm:main}]
 If $\tilde u$ is another minimizer in this homotopy class then equality holds in \eqref{eq:minimality}, \textit{i.e.},  $\tilde u$ is horizontally weakly conformal by Corollary \ref{FS_maps_vs_forms} and $\tilde u^*\omega_{\S^2} \in E^+_{0,1}$, hence $|\tilde u^*\omega_{\S^2}|=4$, cf.\ \eqref{eq:valueIp}. By the invariance of the energy under rotations, using \eqref{SO4_transitive_action}, \eqref{eq:alpha_hopf} and \eqref{eq:basis_for_E_0+}, we can find $\tilde R\in SO(4)$ so that the map $v:=\tilde u\circ \tilde R\in \mathcal{U}_{\mathcal{FS}}$, which is also a minimizer and horizontally weakly conformal, is such that $Q(v)=1$ and $v^*\omega_{\S^2}=h^*\omega_{\S^2}$.  The conclusion then follows from Proposition \ref{prop:weakly_conf_rigid}.
\end{proof}

\section*{Acknowledgments}

AG acknowledges the support of the Royal Society through a Newton International Fellowship. XL was supported by 
the ANR project ANR-22-CE40-0006.  KZ was funded by the Deutsche Forschungsgemeinschaft (DFG, German Research Foundation) – CRC 1720 – 539309657 and also acknowledges previous support by the Hausdorff Center for Mathematics (HCM) under
Germany’s Excellence Strategy -EXC-2047/1-390685813.

The first two authors would also like to thank the Institute for Applied Mathematics of the University of Bonn, as well as the Hausdorff Center for Mathematics in Bonn, for their hospitality during the week 14.04.2025-18.04.2025, when this work was initiated. \BBB XL also thanks the Erwin Schrödinger International Institute for Mathematics and Physics (ESI), where part of this work was conducted during the 2025 thematic program on Free Boundary Problems. \EEE

\section*{Statements}
\textbf{Conflict of interest statement:} The authors declare no conflict of interest.

\textbf{Data availability statement:} No data were used in the preparation of this manuscript.

\bibliographystyle{acm}
\bibliography{ref_hopf}

\end{document}